\documentclass[12pt]{amsart}

\usepackage{makecell}

\usepackage[numbers,sort&compress]{natbib}
\usepackage[colorlinks,
linkcolor=blue,
anchorcolor=green,
citecolor=magenta
]{hyperref}
\usepackage{tikz}
\usepackage{mathrsfs,amssymb,amsfonts}
\usepackage{graphicx}
\usepackage{amsmath,amssymb,latexsym,color}
\usepackage[mathscr]{eucal}
\usepackage[numbers]{natbib}
\usepackage{bm}
\usepackage{subfigure}\usepackage{footnote}
\makeatletter

\numberwithin{equation}{section}

\newtheorem{theorem}{Theorem}[section]
\newtheorem{claim}{Claim}[]
\newtheorem{definition}[theorem]{Definition}

\newtheorem{lemma}[theorem]{Lemma}
\newtheorem{problem}[theorem]{Problem}
\newtheorem{remark}[theorem]{Remark}

\newtheorem{corollary}[theorem]{Corollary}

\begin{document}
	\title[Disjoint chorded cycles]{\bf 
		Disjoint chorded cycles in a $2$-connected graph\textsuperscript{*}}
	\thanks{*Supported by the National Natural Science Foundation of China (12471328, 12331013) and the Fundamental Research Funds for the Central Universities.}
	
	\author[Lu]{Zai Ping Lu} 
	\address{
		Z.P. Lu\\
		Center for Combinatorics, LPMC,
		Nankai University,
		Tianjin 300071, China
	}
	\email{lu@nankai.edu.cn}
	\author[Xue]{Shu Dan Xue\textsuperscript{**}}\thanks{**Corresponding author. }
	\address{
		S.D. Xue\\
		Center for Combinatorics, LPMC,
		Nankai University,
		Tianjin 300071, China
	}
	\email{1120220002@mail.nankai.edu.cn}

	\begin{abstract} 
		A  chorded cycle in a graph $G$
		is a cycle containing an edge of $G$ that joins two nonconsecutive vertices of the cycle.
		In 2010, Gao and Qiao  independently proved that
		a graph of order at least $4s$, in which the neighborhood union of any two nonadjacent vertices has at least $4s+1$ vertices, contains $s$ vertex-disjoint chorded cycles.
		In 2022, Gould raised a problem that asks whether increasing connectivity would improve the neighborhood union condition.
    	In this paper, we solve the problem for $2$-connected graphs by proving
    	that a $2$-connected graph of order at least $4s$, in which the
    	neighborhood union of any two nonadjacent vertices has at least $4s$
	    vertices, contains $s$ vertex-disjoint chorded cycles. Moreover, the
    	neighborhood-union bound $4s$ is sharp.

		\vskip 10pt
		
		
	\end{abstract}
	\maketitle
	\date\today

	
	\parskip=2pt
	
	\section{Introduction}\label{introduction}

	In this paper, all graphs are assumed to be finite and simple.
	
	Let $G$ be a graph with vertex set $V(G)$ and edge set $E(G)$. The neighborhood and degree of a vertex $u$ in $G$ are denoted by $N_G(u)$ and $deg_G(u)$, respectively. 
	For a subset $S\subseteq V(G)$ and
	an integer $m\geqslant 1$, put
	\[
	\begin{split}
		N_G(S)&=\left\{u\in V(G)\,:\, \{u,v\}\in E(G)\mbox{ for some }v\in S\right\},\\
		\sigma_m(G)&=\min\left\{\sum_{u\in S}deg_G(u)\,:\,  S \mbox{ is an independent set of size } m\right\},\\
		\delta_m(G)&=\min\left\{|N_G(S)|\,:\,  S \mbox{ is an independent set of size } m\right\}.
	\end{split}
	\]
	Note that $\sigma_1(G)=\delta_1(G)$ is just the minimum degree $\delta(G)$ of $G$.

	A {\em chord} of a cycle $C$ in a graph $G$ is an edge in $E(G)\setminus E(C)$ both of whose ends lie on 
	$C$. A {\em chorded cycle} is a cycle that has at least one chord.
	Exploring conditions  on $\delta_m(G)$, $\sigma_m(G)$ or $|E(G)|$
	that guarantee a graph $G$ has $s$ vertex-disjoint chorded cycles is a fascinating and challenging problem. Table \ref{tab:1} summarizes some of the latest   results on the existence of $s$ vertex-disjoint chorded cycles in a graph.
	
	
	The following was shown independently in \cite{YGJ} and \cite{SNQ}.
	
	\begin{theorem}\label{T2}
		Let $G$ be a graph of order at least $4s$, where $s \geqslant 1$. If $\delta_{2}(G) \geqslant 4s+1$, then $G$ contains $s$ vertex-disjoint chorded cycles.
	\end{theorem}
	
	Gould \cite{RJG1} raised the question of whether increasing connectivity would improve the outcome.
	
	\begin{problem}\label{P1}
		Can $\delta_{2}(G)$ of Theorem {\rm \ref{T2}} be decreased if the graph $G$ is  $k$-connected for some $k \geqslant 2$?
	\end{problem}

	\begin{table}
		\caption{Previously known results}
		\centering
		\begin{tabular*}{14cm}{cccccccc}
			\hline
			Condition	&  ~~$\delta(G)$ & ~~~$\sigma_{2}(G)$ & ~~~$\sigma_{3}(G)$  & ~~~$\sigma_{4}(G)$  &~~~$\sigma_{m}(G)$ & ~~~$\delta_{2}(G)$  \\
			\hline
			\\
			Lower bound	& ~~$3s$ & ~~~$6s-1$ & ~~~$9s-2$  & ~~~$12s-3$  & ~~~$3sm-m+1$ & ~~~$4s+1$   \\
			\\
			Reference &  ~~\cite{DF} &~~~ \cite{CFGL} & ~~~\cite{RKA}  &  ~~~\cite{RKA1}  & ~~~\cite{BRK} & ~~~ \cite{YGJ,SNQ} \\
			\hline
			
		\end{tabular*}
	\end{table}\label{tab:1}

	We consider here the case that $k=2$, and give the following theorem. 
	
	\begin{theorem}\label{T4}
		Let  $G$ be a $2$-connected graph of order at least $4s$, where $s \geqslant 1$. Suppose that either $\delta_{2}(G) \geqslant 4s$ or $G$ is a complete graph. Then $G$ contains $s$ vertex-disjoint chorded cycles.
	\end{theorem}
	


			

	The following remark says that the bound $4s$ for $\delta_2(G)$ in Theorem \ref{T4} is sharp.

	\begin{remark}	
		\normalfont 
    For every integer $s\geqslant 1$, let $G_s$ be the graph obtained from the complete graph $K_{4s-1}$ by subdividing a given edge $ab$ twice, so that the edge $ab$ is replaced by the path $axyb$. Then $G_s$ is a $2$-connected graph of order $4s+1$.

      Given any pair of nonadjacent vertices $u,v\in V(G_s)$, every vertex in $V(G_s)\setminus\{u,v\}$ is adjacent to at least one of $u$ and $v$, and thus
		\[
		N_{G_s}(u)\cup N_{G_s}(v)
		=V(G_s)\setminus\{u,v\}.
		\]
		Hence,
		\[
		\delta_2(G_s)=|V(G_s)|-2=4s-1.
		\]
		
		Suppose that $\mathcal{F}$ is a collection
		of $s$ vertex-disjoint chorded cycles in $G_s$. Let $U$ be the vertex set of $K_{4s-1}$.
		Since every chorded cycle has length at least $4$ and $|U|=4s-1$, some
		$C\in\mathcal{F}$ must contain $x$ or $y$. Since
		$\deg_{G_s}(x)=\deg_{G_s}(y)=2$,
		the cycle $C$ must contain the path $axyb$. Therefore, $C$ is the union of the path $axyb$ and a path $P$ from $b$ to $a$ whose internal vertices lie in $U\setminus\{a,b\}$.
		
		Since $ab\notin E(G_s)$, the path $P$ has at least one internal vertex. If $P$ has exactly one internal vertex $z$,
		then $C=axybza$, and so $C$ is not a chorded cycle, a contradiction.

		This implies that the path $P$ has at least two
		internal vertices, and hence $C$ contains at least four vertices in
		$U$. Moreover, $C$ contains
		both $x$ and $y$. Therefore, each of the remaining $s-1$ chorded
		cycles in $\mathcal{F}$ has all its vertices in $U$ and contains at
		least four further vertices of $U$. Consequently, the cycles in
		$\mathcal{F}$ use at least $4+4(s-1)=4s$ vertices of $U$, contradicting $|U|=4s-1$. Therefore, $G_s$ contains no $s$ vertex-disjoint chorded cycles. 
		Hence the bound $4s$ in Theorem~\ref{T4} is sharp. 

\end{remark}
	\vskip 20pt

	\section{Chorded cycles in a  graph}
	\label{sec:pre}
	In this section, we make some preparations   for the proof
	of Theorem \ref{T4} by collecting several known results and proving some technical lemmas, which involve either constructing a chorded cycle or proving its existence.

	We first introduce the notation used in this and the following sections.
	
	Let $G$ be a graph. 
	An edge $\{u,v\}$ of $G$ is always viewed as a path of length $1$ and written as $uv$. A path or cycle  of $G$ of length $\ell$ is always written as a sequence $u_1u_2\cdots u_{\ell+1}$ of   vertices with $u_{i}u_{i+1}\in E(G)$ for all $1\leqslant i\leqslant \ell$ and,  in the cycle case, $u_{\ell+1}=u_1$. For a subset $S\subseteq V(G)$, denote by $\langle S\rangle$ the subgraph of $G$ induced by $S$, and put $G-S=\langle V(G)\setminus S\rangle$ (if $S\ne V(G)$). When $S$ is a singleton, say $S=\{u\}$, we write $\langle S\rangle$  and $G-S$ simply as $u$ and $G-u$, respectively. In addition, for a subgraph $H$ of $G$ with $V(H)\ne V(G)$, put $G-H=\langle V(G)\setminus V(H)\rangle$.
	
	Let $H$ be a subgraph of   $G$. If $u\in V(G)$, then $N_H(u)$ denotes the set of neighbors of $u$ that lie in $H$ , that is, $N_H(u)=N_G(u)\cap V(H)$, and put $d_H(u)=|N_H(u)|$.
	(Note that even when  $u\in V(H)$, the value of $d_H(u)$ may be larger than the degree $deg_H(u)$ of $u$ in $H$.) 
		For two vertices $u,v\in V(G)$, define
	$d_H(u,v)=|N_H(u)\cup N_H(v)|$.
		Thus, $d_H(u,v)$ denotes the size of the union of the neighborhoods of $u$ and $v$ in $H$, not the distance between $u$ and $v$.
	Similarly, for $S\subseteq V(G)$, put $N_H(S)=N_G(S)\cap V(H)$, and put $d_H(S)=|N_H(S)|$.
	If $X,Y\subseteq V(G)$ or $X$ and $Y$ are subgraphs of $G$ then $E_H(X,Y)$, written as $E(X,Y)$ when $H=G$, denotes the  set of edges of $H$ connecting a vertex in $X$ and a vertex in $Y$. 
	If $X$ and $Y$ are subgraphs of $G$, then $X\sqcup_H Y$, written simply as $X\sqcup Y$ when $H=G$, denotes the subgraph with vertex set $V(X)\cup V(Y)$ and edge set $E(X)\cup E(Y)\cup E_H(X,Y)$.

	Let $H$ be either a path or a cycle of a graph $G$. For vertices $u,v \in V(H)$, the notation  $H[u,v]$  stands for a path of $H$  that connects $u$ and $v$. Clearly,  $H[u,v]$ is uniquely determined when  $H$ is a path, and there are two choices of $H[u,v]$ when  $H$ is a cycle (of positive length). 
  In the latter case, we define $H[u,v]$ as follows. Label  the vertices  of the cycle as $H=u_1 u_2\cdots u_i\cdots  u_\ell u_1$, if $i\leqslant j$ then  $H[u_i,u_j]$  stands for the path $u_iu_{i+1}\cdots u_{j-1}u_j$, while $H[u_j,u_i]$
	is the path $u_ju_{j+1}\cdots u_\ell u_1\cdots u_{i-1}u_i$.

	\vskip 5pt
	
	The following lemma presents some sufficient conditions for the existence of a chorded cycle,
	which can be deduced from \cite[Lemmas 3.4, 3.5, and 3.7]{RKA}.

	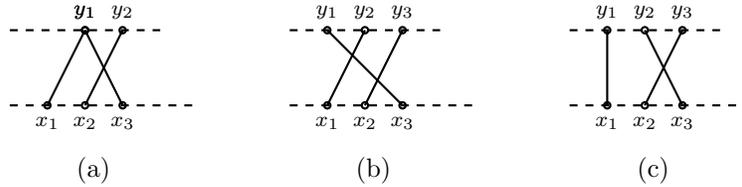
\begin{figure}[ht]
		\centering
		\subfigure[]{
			\begin{tikzpicture}[style=thick]
				\draw[dashed] (0,0)  -- (2.5,0);
				\draw[dashed] (0,1)   -- (2,1);
				\draw (0.5,0) circle (1.2pt) node [below]{\tiny $x_1$}-- (1,1)circle (1.2pt)node [above]{\tiny $y_1$};
				\draw (1.5,0)node [below]{\tiny $x_3$} circle (1.2pt) -- (1,1)circle (1.2pt)node [above]{\tiny $y_1$};
				\draw (1,0) circle (1.2pt) node [below]{\tiny $x_2$}-- (1.5,1)circle (1.2pt)node [above]{\tiny $y_2$};
		\end{tikzpicture}}\hskip 30pt 
		\subfigure[]{
			\begin{tikzpicture}[style=thick]
				\draw[dashed] (0,0)  -- (2.5,0);
				\draw[dashed] (0,1)   -- (2,1);
				\draw (0.5,0) circle (1.2pt) node [below]{\tiny $x_1$}-- (1,1)circle (1.2pt)node [above]{\tiny $y_2$};
				\draw (1.5,0)node [below]{\tiny $x_3$} circle (1.2pt) -- (0.5,1)circle (1.2pt)node [above]{\tiny $y_1$};
				\draw (1,0) circle (1.2pt) node [below]{\tiny $x_2$}-- (1.5,1)circle (1.2pt)node [above]{\tiny $y_3$};
		\end{tikzpicture}}\hskip 30pt 
		\subfigure[]{
			\begin{tikzpicture}[style=thick]
				\draw[dashed] (0,0)  -- (2.5,0);
				\draw[dashed] (0,1)   -- (2,1);
				\draw (0.5,0) circle (1.2pt) node [below]{\tiny $x_1$} -- (0.5,1)circle (1.2pt)node [above]{\tiny $y_1$};
				\draw (1.5,0) circle (1.2pt) node [below]{\tiny $x_3$} -- (1,1)circle (1.2pt)node [above]{\tiny $y_2$};
				\draw (1,0) circle (1.2pt)  node [below]{\tiny $x_2$}-- (1.5,1)circle (1.2pt)node [above]{\tiny $y_3$};
		\end{tikzpicture}} 
		\caption{Exceptions for  $|E_H(V_1,V_2)|=3$}\label{fig-1}
	\end{figure}

	\begin{figure}[ht]
		\centering
		\subfigure[]{
			\begin{tikzpicture}[style=thick]
				\draw[dashed] (0,0)  -- (2.5,0);
				\draw[dashed] (0,1)   -- (2,1);
				\draw (0.5,0) circle (1.2pt) node [below]{\tiny $x_1$}-- (1,1)circle (1.2pt)node [above]{\tiny $y_2$};
				\draw (1.5,0) circle (1.2pt) node [below]{\tiny $x_3$}-- (1,1)circle (1.2pt)node [above]{\tiny $y_2$};
				\draw (1,0) circle (1.2pt) node [below]{\tiny $x_2$}-- (1.5,1)circle (1.2pt)node [above]{\tiny $y_3$};
				\draw (1,0) circle (1.2pt) node [below]{\tiny $x_2$}-- (0.5,1)circle (1.2pt)node [above]{\tiny $y_1$};
		\end{tikzpicture}}\hskip 30pt
		\subfigure[]{
			\begin{tikzpicture}[style=thick]
				\draw[dashed] (0,0)  -- (2.5,0);
				\draw[dashed] (0,1)   -- (2,1);
				\draw (0.3,0) circle (1.2pt)node [below]{\tiny $x_1$} -- (1,1)circle (1.2pt)node [above]{\tiny $y_2$};
				\draw (1.7,0) circle (1.2pt) node [below]{\tiny $x_4$}-- (1,1)circle (1.2pt)node [above]{\tiny $y_2$};
				\draw (0.8,0) circle (1.2pt) node [below]{\tiny $x_2$}-- (1.5,1)circle (1.2pt)node [above]{\tiny $y_3$};
				\draw (1.2,0) circle (1.2pt) node [below]{\tiny $x_3$}-- (0.5,1)circle (1.2pt)node [above]{\tiny $y_1$};
		\end{tikzpicture}}\hskip 30pt
		\subfigure[]{
			\begin{tikzpicture}[style=thick]
				\draw[dashed] (0,0)  -- (2.5,0);
				\draw[dashed] (0,1)   -- (2,1);
				\draw (0.3,0) circle (1.2pt)node [below]{\tiny $x_1$} -- (1,1)circle (1.2pt)node [above]{\tiny $y_2$};
				\draw (1.7,0) circle (1.2pt)node [below]{\tiny $x_4$} -- (1,1)circle (1.2pt)node [above]{\tiny $y_2$};
				\draw (1.2,0) circle (1.2pt) node [below]{\tiny $x_3$}-- (1.5,1)circle (1.2pt)node [above]{\tiny $y_3$};
				\draw (0.8,0) circle (1.2pt) node [below]{\tiny $x_2$}-- (0.5,1)circle (1.2pt)node [above]{\tiny $y_1$};
		\end{tikzpicture}}

		\subfigure[]{
			\begin{tikzpicture}[style=thick]
				\draw[dashed] (0,0)  -- (2.5,0);
				\draw[dashed] (0,1)   -- (2,1);
				\draw (0.8,0) circle (1.2pt) node [below]{\tiny $x_2$} -- (0.3,1)circle (1.2pt)node [above]{\tiny $y_1$};
				\draw (0.3,0) circle (1.2pt) node [below]{\tiny $x_1$} -- (0.8,1)circle (1.2pt)node [above]{\tiny $y_2$};
				\draw (1.7,0) circle (1.2pt) node [below]{\tiny $x_4$} -- (1.2,1)circle (1.2pt)node [above]{\tiny $y_3$};
				\draw (1.2,0) circle (1.2pt) node [below]{\tiny $x_3$} -- (1.7,1)circle (1.2pt)node [above]{\tiny $y_4$};
		\end{tikzpicture}}\hskip 30pt
		\subfigure[]{
			\begin{tikzpicture}[style=thick]
				\draw[dashed] (0,0)  -- (2.5,0);
				\draw[dashed] (0,1)   -- (2,1);
				\draw (1.3,0) circle (1.2pt)node [below]{\tiny $x_3$}  -- (0.3,1)circle (1.2pt)node [above]{\tiny $y_1$};
				\draw (0.3,0) circle (1.2pt) node [below]{\tiny $x_1$} -- (0.8,1)circle (1.2pt)node [above]{\tiny $y_2$};
				\draw (1.8,0) circle (1.2pt)node [below]{\tiny $x_4$}  -- (1.3,1)circle (1.2pt)node [above]{\tiny $y_3$};
				\draw (0.8,0) circle (1.2pt) node [below]{\tiny $x_2$} -- (1.8,1)circle (1.2pt)node [above]{\tiny $y_4$};
		\end{tikzpicture}}\hskip 30pt
		\subfigure[]{
			\begin{tikzpicture}[style=thick]
				\draw[dashed] (0,0)  -- (2.5,0);
				\draw[dashed] (0,1)   -- (2,1);
				\draw (1.3,0) circle (1.2pt) node [below]{\tiny $x_3$}-- (0.3,1)circle (1.2pt)node [above]{\tiny $y_1$};
				\draw (0.3,0) circle (1.2pt) node [below]{\tiny $x_1$}-- (1.3,1)circle (1.2pt)node [above]{\tiny $y_3$};
				\draw (1.8,0) circle (1.2pt)node [below]{\tiny $x_4$} -- (0.8,1)circle (1.2pt)node [above]{\tiny $y_2$};
				\draw (0.8,0) circle (1.2pt)node [below]{\tiny $x_2$} -- (1.8,1)circle (1.2pt)node [above]{\tiny $y_4$};
		\end{tikzpicture}}
		\caption{Exceptions for $|E_H(V_1,V_2)|=4$}\label{fig-2}
	\end{figure}
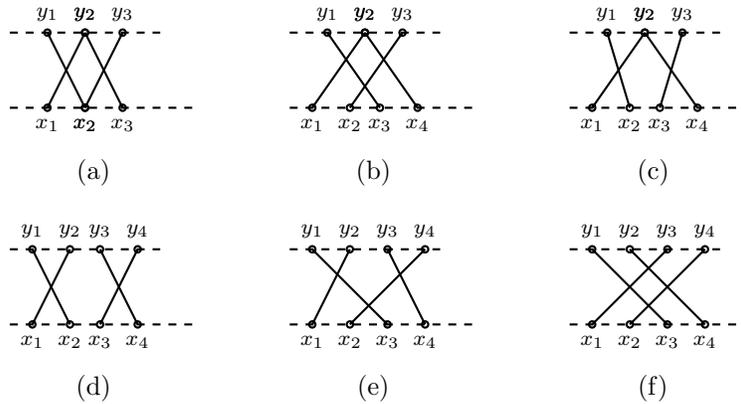

	\begin{lemma}\label{vertex-disjoint-paths}
		Let  $H$ be a graph with vertex set partitioned into two nonempty sets $V_1$ and $V_2$ such that both $\langle V_1\rangle$ and $\langle V_2\rangle$ are paths.  Then  $H$ contains no chorded cycles if and only if either $|E_H(V_1,V_2)|\leqslant 2$ or $H$ is  isomorphic to one of the  graphs illustrated in Figures {\rm \ref{fig-1}} and {\rm \ref{fig-2}}. In particular, if $|E_H(V_1,V_2)|\geqslant 5$ then $H$ contains a chorded cycle.
	\end{lemma}
	\proof 
	Clearly, if either $|E_H(V_1,V_2)|\leqslant 2$ or $H$ is isomorphic to one of the graphs illustrated in Figures {\rm \ref{fig-1}} and {\rm \ref{fig-2}}, then $H$ contains no chorded cycles. Next we suppose  that $|E_H(V_1,V_2)|\geqslant 3$,  and show that either $H$ contains a chorded cycle or $H$ is   isomorphic to one of the   graphs illustrated in Figures {\rm \ref{fig-1}} and {\rm \ref{fig-2}}.

\begin{figure}[ht]
		\centering
		\subfigure[]{
			\begin{tikzpicture}[style=thick]
				\draw[dashed] (0,0)  -- (2.0,0);
				\draw[dashed] (0,1)   -- (2,1);
				\draw (0.5,0) circle (1.2pt)  -- (1,1)circle (1.2pt);
				\draw (1,0)circle (1.2pt)  -- (1,1)circle (1.2pt);
				\draw (1.5,0) circle (1.2pt) -- (1,1)circle (1.2pt);
		\end{tikzpicture}}\hskip 20pt
		\subfigure[]{
			\begin{tikzpicture}[style=thick]
				\draw[dashed] (0,0)  -- (2.0,0);
				\draw[dashed] (0,1)   -- (2,1);
				\draw (0.5,0) circle (1.2pt)  -- (1,1)circle (1.2pt);
				\draw (1,0)circle (1.2pt)  -- (1.5,1)circle (1.2pt);
				\draw (1,0) circle (1.2pt) -- (1,1)circle (1.2pt);
		\end{tikzpicture}}\hskip 20pt
		\subfigure[]{
			\begin{tikzpicture}[style=thick]
				\draw[dashed] (0,0)  -- (2.0,0);
				\draw[dashed] (0,1)   -- (2,1);
				\draw (0.5,0) circle (1.2pt)  -- (1,1)circle (1.2pt);
				\draw (1.5,0)circle (1.2pt)  -- (1.5,1)circle (1.2pt);
				\draw (1,0) circle (1.2pt) -- (1,1)circle (1.2pt);
		\end{tikzpicture}}

		\subfigure[]{
			\begin{tikzpicture}[style=thick]
				\draw[dashed] (0,0)  -- (2.0,0);
				\draw[dashed] (0,1)   -- (2,1);
				\draw (0.5,0) circle (1.2pt)  -- (1.5,1)circle (1.2pt);
				\draw (1.5,0)circle (1.2pt)  -- (.5,1)circle (1.2pt);
				\draw (1,0) circle (1.2pt) -- (1.5,1)circle (1.2pt);
		\end{tikzpicture}}\hskip 20pt
		\subfigure[]{
			\begin{tikzpicture}[style=thick]
				\draw[dashed] (0,0)  -- (2.0,0);
				\draw[dashed] (0,1)   -- (2,1);
				\draw (0.4,0) circle (1.2pt)  -- (1.8,1)circle (1.2pt);
				\draw (1.8,0)circle (1.2pt)  -- (.4,1)circle (1.2pt);
				\draw (1,0) circle (1.2pt) -- (1,1)circle (1.2pt);
		\end{tikzpicture}}\hskip 20pt
		\subfigure[]{
			\begin{tikzpicture}[style=thick]
				\draw[dashed] (0,0)  -- (2.0,0);
				\draw[dashed] (0,1)   -- (2,1);
				\draw (0.5,0) circle (1.2pt)  -- (0.5,1)circle (1.2pt);
				\draw (1.5,0)circle (1.2pt)  -- (1.5,1)circle (1.2pt);
				\draw (1,0) circle (1.2pt) -- (1,1)circle (1.2pt);
		\end{tikzpicture}} 
		\caption{Chorded cycles in $H$}\label{fig-3}
	\end{figure}
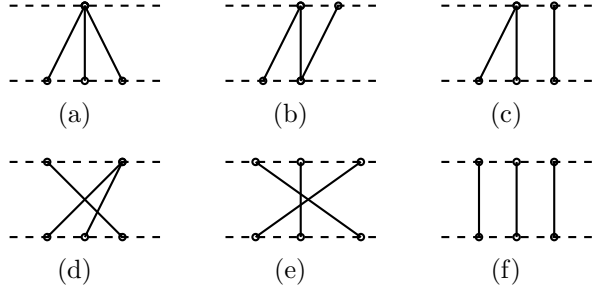

	If $|E_H(V_1, V_2)|\in\{3,4\}$, then it is straightforward to check that, except the graphs illustrated in Figures {\rm \ref{fig-1}} and {\rm \ref{fig-2}}, $H$ has a subgraph isomorphic to one of the graphs illustrated in Figure \ref{fig-3}, where each graph contains a chorded cycle. 
	
	Thus we suppose further that $|E_H(V_1,V_2)|\geqslant 5$. Let $F\subseteq E_H(V_1,V_2)$, and denote by $H(F)$ the subgraph  of $H$ obtained from $\langle V_1\rangle\cup \langle V_2\rangle$ by adding the edges in $F$.
	Choose $F$ with $|F|=4$. If $H(F)$ contains a chorded cycle, then so does $H$. Suppose that $H(F)$ contains no chorded cycles. Then $H(F)$ is isomorphic to one of the graphs illustrated in Figure  {\rm \ref{fig-2}}.  Pick $uv\in E_H(V_1,V_2)\setminus F$, and let $E=F\cup\{uv\}$. It is straightforward to check, by analyzing the locations of $u$ and $v$,  that $H(E)$ has a subgraph isomorphic to one of the graphs illustrated in Figure  {\rm \ref{fig-3}}. Then $H(E)$ contains a chorded cycle, and hence so does $H$. This completes the proof.
	\qed
	
	\vskip 10pt
	
	Clearly, if a subgraph of a graph $G$ contains a chorded cycle, then so does $G$. This leads to the following simple observations concerning the degree of special vertices when a graph  does not contain a chorded cycle, see also \cite{RKA,RKA1}.

	\begin{lemma}
		\label{v2}
		Let $H$ be a graph without chorded cycles. Suppose that $H$ contains a path $P=u_{1}u_{2} \cdots u_{p}$, where $p\geqslant  3$. 
		\begin{itemize}
			\item[\rm (1)] If $u_{1}u_{i} \in E(H)$ with $i \geqslant 3$, then $d_{P}(u_{j}) \leqslant 3$ for all $j \leqslant i-1$,  and $d_{P}(u_{i-1})=2$.
			\item[\rm (2)] If $u_{p}u_{i} \in E(H)$ with $i \leqslant p-2$, then $d_{P}(u_{j}) \leqslant 3$ for all $j \geqslant i+1$,  and $d_{P}(u_{i+1})=2$.
		\end{itemize}
	\end{lemma}
	
	\begin{lemma}
		\label{degree-2}
		Let $H$ be a connected graph without chorded cycles and  Hamiltonian paths. Suppose that
		$P_{1}=u_{1}u_{2} \cdots u_{p}$ is a longest path in $H$ with $p\geqslant 3$, and   $P_{2}=v_{1}v_{2} \cdots v_{q}$ is a longest path in $H - P_{1}$  with $q\geqslant 1$.  
		Then the following statements hold.
		\begin{itemize}
			\item [\rm (1)] If $i \in \{1,p\}$ then $d_{H-P_{1}}(u_{i})=0$.
			\item [\rm (2)] If $i \in \{1,p\}$ then $d_{H}(u_{i})=d_{P_{1}}(u_{i}) \leqslant 2$.
			\item [\rm (3)] If $j \in \{1,q\}$ then $d_{H- (P_{1} \cup P_{2})}(v_{j})= 0$.
			\item [\rm (4)] If $j \in \{1,q\}$ then $d_{P_{2}}(v_{j}) \leqslant 2$.
			\item [\rm (5)] If $i \in \{1,2\}$ and $w\in V(H) \setminus V(P_{i})$  then $d_{P_{i}}(w) \leqslant 2$.
			\item [\rm (6)] If $q \geqslant 2$ then $d_{P_{1}}(v_{1})+d_{P_{1}}(v_{q}) \leqslant 3$.
		\end{itemize}
	\end{lemma}

	\begin{lemma}\label{no-Hamiltonian-path}
		Let $H$ be a connected graph without chorded cycles and  Hamiltonian paths. 
		Suppose that  $P_{1}=u_{1}u_{2} \cdots u_{p}$ is a longest path in $H$ with $p\geqslant 3$, and   $P_{2}=v_{1}v_{2} \cdots v_{q}$ is a longest path in $H - P_{1}$ with $q\geqslant 1$.  Suppose that $|V(H)|\geqslant 4$ and
		$d_{P_{1}}(v_{1}) \leqslant d_{P_{1}}(v_{q})$.
		Then 
		\begin{itemize}
			\item [\rm (1)] $q\leqslant 2$, and $V(H)= V(P_{1}) \cup  V(P_{2})$; or
			\item [\rm (2)] $q \geqslant 3$, and $d_{H}(v_{1})=1$; or
			\item [\rm (3)] there exists   $w\in V(H-(P_{1}\cup  v_{1}))$ such that $d_{H}(w) \leqslant 2$, $u_{1}w, u_{p}w \notin E(H)$ and $w$ is not a cut-vertex of $H$. 
		\end{itemize}
	\end{lemma}
	\proof 
	Since $H$ is connected and contains no Hamiltonian paths, there exist $u_i\in V(P_1)$ and $v\in V(H-P_1)$ with $u_iv\in E(H)$.  
	If $u_{1}u_{p} \in E(H)$, then there is a longer path $vu_iu_{i-1}\cdots u_1u_p\cdots u_{i+1}$ than $P_{1}$, which contradicts the choice of $P_1$. 
	Thus $u_{1}u_{p} \notin E(H)$. 
	Also, by the choice of $P_1$, we have $u_1v, u_pv\not\in E(H)$ for any  $v\in V(H-P_1)$.  
	We next discuss two cases according to $q\leqslant 2$ and $q\geqslant 3$, respectively.
	
	{\it Case}~1. Suppose that $q \leqslant 2$. If $V(H)= V(P_{1}) \cup  V(P_{2})$ then (1) of the lemma occurs.
	Next we suppose that $V(H)\ne V(P_{1}) \cup  V(P_{2})$.
	Put $K=H- (P_{1} \cup P_{2})$.  Then, since $q=1$ or $2$, we have
	$d_{K}(v) =0$ for all $v \in V(P_{2})$, see (3) of Lemma \ref{degree-2}.
	Since $H$ is connected, this implies that  $N_{K}(u) \neq \emptyset$ for some $u\in V(P_1)$. 
	
	Pick $w_{1} \in N_{K}(u)$. Then $u_{1}w_1,u_{p}w_1 \notin E(H)$,   $d_{P_2}(w_1)=0$, and $d_H(w_1)=d_{K}(w_{1})+d_{P_{1}}(w_{1})$.  Recall from (5) of Lemma \ref{degree-2} that $d_{P_1}(w_1) \leqslant 2$.
	If $d_{K}(w_{1})=0$ then $d_H(w_1)=d_{P_1}(w_1) \leqslant 2$, and so $w_1$ is not a cut-vertex of $H$. Taking  $w=w_{1}$, (3) of the lemma occurs.
	Thus $d_{K}(w_{1}) \geqslant 1$, and we put $N_{K}(w_{1}) =\{w_2\}$.
	If $d_{K}(w_{1}) \geqslant 2$ or $d_{K}(w_{2}) \geqslant 2$, then $K$ has a path with at least three vertices, which contradicts the choice of $P_{2}$. This says that $d_{K}(w_{i}) =1$, and so $d_{P_{1}}(w_{i})=d_{H}(w_{i})-d_{K}(w_{i})=d_{H}(w_{i})-1$ for each $i \in \{1,2\}$. 
	Without loss of generality, we assume that $d_{H}(w_{1}) \geqslant d_{H}(w_{2})$, that is, $d_{P_1}(w_{1}) \geqslant d_{P_1}(w_{2})$.
	If $d_{H}(w_{2}) \leqslant 2$ then $u_{1}w_2,u_{p}w_2 \notin E(H)$, $w_2$ is not a cut-vertex of $H$ and (3) of the lemma  occurs by taking $w=w_{2}$.
	Thus $d_{H}(w_{1}) \geqslant d_{H}(w_{2}) \geqslant 3$, and further $d_{P_{1}}(w_{i})=d_{H}(w_{i})-1 \geqslant 2$ for each $i \in \{1,2\}$.  Considering the subgraph $P_1 \sqcup w_1w_2$, it follows from  Lemma \ref{vertex-disjoint-paths} that $H$ contains a chorded cycle, a contradiction.

	{\it Case}~2. Suppose that  $q \geqslant 3$. 
	Since $d_{P_{1}}(v_{1}) \leqslant d_{P_{1}}(v_{q})$, by (5) and (6) of Lemma \ref{degree-2}, $d_{P_{1}}(v_{1}) \leqslant 1$ and $d_{P_{1}}(v_{q}) \leqslant 2$.  
	Suppose that $d_{P_{1}}(v_{q}) =0$, and so  $d_{P_{1}}(v_{1})=0$. Then $d_{H}(v_{q})=d_{P_{2}}(v_{q})$, and so $v_q$ is not a cut-vertex of $H$.
	According to (3) and (4) of Lemma \ref{degree-2}, we deduce that  $d_{H}(v_{j}) \leqslant 2$, where $j\in \{1,q\}$. 
	Recalling that  $u_1v_q, u_pv_q\not\in E(H)$,  (3) of the lemma occurs by taking $w=v_q$. Thus, in the following, we let $1\leqslant d_{P_{1}}(v_{q}) \leqslant  2$  and 
	$d_{P_{1}}(v_{1}) \leqslant 1$.

	{\it Subcase} 2.1. Suppose that   $d_{P_{1}}(v_{1}) =1$.
	If $d_{P_{2}}(v_{1})\geqslant 2$ or $d_{P_{2}}(v_{q})\geqslant 2$ or  $d_{P_{2}}(v_{j})\geqslant 3$ for some $2\leqslant j\leqslant q-1$, then there exists a chorded cycle  in $P_1\sqcup \langle V(P_2)\rangle$ with a chord adjacent to $v_{1}$ or $v_{q}$ or $v_j$, respectively,  a contradiction.  This forces that $d_{P_{2}}(v_{1})=d_{P_{2}}(v_{q})=1$, and $d_{P_{2}}(v_{j})=2$ for all $2\leqslant j\leqslant q-1$. In particular, $\langle V(P_2)\rangle =P_2$.
	
	If  $d_{P_{1}}(v_{q})=1$ then, combining (3) of Lemma \ref{degree-2},  $d_{H}(v_{q})=2$, $v_q$ is not a cut-vertex of $H$ and so (3) of the lemma occurs by taking $w=v_q$. 
	If $d_H(v_{q-1})=d_{P_2}(v_{q-1})= 2$,  then $v_{q-1}$ is not a cut-vertex of $H$ and  so (3) of the lemma occurs by taking $w=v_{q-1}$.  
	Thus, we suppose next that $d_{P_{1}}(v_{q})=2$ and $d_H(v_{q-1})\geqslant 3$.
	In addition, since $d_{P_{1}}(v_{1}) =1$,  we have $d_{H}(v_{1})=2$ by (3) of Lemma \ref{degree-2}.
	
	Considering the subgraph  $P_1\sqcup P_2$, since $H$ contains no chorded cycles, it follows from Lemma \ref{vertex-disjoint-paths} that $|E_H(P_1,P_2)|=3$, and the subgraph $P_1\sqcup P_2$ is described as in (a) of Figure \ref{fig-1} with $v_q=y_1$, $v_1=y_2$, $N_{P_1}(v_q)=\{x_1,x_3\}$ and $N_{P_1}(v_1)=\{x_2\}$.
	In particular, $N_{P_1}(v_{q-1})=\emptyset$.
	Recalling that $d_{P_{2}}(v_{q-1})=2$ and $d_H(v_{q-1})\geqslant 3$,  there exists $w_1\in V(H-(P_1\cup P_2))$ with $w_1v_{q-1}\in E(H)$.
	If  $d_{P_1}(w_1)\geqslant 2$  or $d_{P_2}(w_1)\geqslant 3$
	then, by Lemma \ref{vertex-disjoint-paths}, either $P_1\sqcup v_qv_{q-1}w_1$ or $P_2\sqcup w_1$ contains a chorded cycle, a contradiction. 
	Thus $d_{P_1}(w_1)\leqslant 1$ and $d_{P_2}(w_1)\leqslant 2$.
	Moreover, $N_H(w_1)\subseteq V(P_1)\cup V(P_2)$; otherwise, $H-P_1$ has a path with at least $q+1$ vertices, which contradicts the choice of $P_2$. Then $d_H(w_1)=d_{P_1}(w_1)+d_{P_2}(w_1)\leqslant 3$.
	
	Suppose that $d_H(w_1)\geqslant 3$. This forces that $d_H(w_1)= 3$,  $d_{P_1}(w_1)=1$ and $d_{P_2}(w_1)=2$. By (3) of Lemma \ref{degree-2}, we have $v_1, v_q\not\in N_{P_2}(w_1)$.
	Put $N_{P_1}(v_1)=\{u_i\}$, $N_{P_1}(w_1)=\{u_j\}$ and $N_{P_2}(w_1)=\{v_k,v_{q-1}\}$, 
	where $2 \leqslant  k\leqslant  q-2$. Let $P=P_1[u_i,u_j]$ and $Q=v_1\cdots v_{q-1}w_1$. 
	Then the subgraph $P\sqcup Q$ contains a chorded cycle with chord $w_{1}v_{k}$, a contradiction.
	Therefore, $d_H(w_1)\leqslant 2$. Clearly, $H-w_1$ is connected.
	By (1) of Lemma \ref{degree-2}, we have $u_1w_1, u_pw_1\not\in E(H)$, and 
	then (3) of the lemma occurs by taking $w=w_1$.

	{\it Subcase}~2.2. Suppose that  $d_{P_{1}}(v_{1}) =0$.
	By (3) and (4) of Lemma \ref{degree-2}, $d_{H}(v_{1})=d_{P_{2}}(v_{1})\leqslant 2$. 
	If $d_{H}(v_{1})=d_{P_{2}}(v_{1})=1$ then (2) of the lemma follows. We next let $d_{H}(v_{1})=d_{P_{2}}(v_{1})=2$, and put $N_{P_2}(v_1)=\{v_2,v_k\}$. 
	
	Clearly, $3\leqslant k\leqslant q$. 
	If $d_{P_1}(v_{k-1})\ne 0$ then, recalling that $d_{P_{1}}(v_{q})\geqslant 1$, it is easily shown that $P_1\sqcup \langle V(P_2)\rangle$ contains a 
	chorded cycle with chord $v_{k-1}v_k$, which gives rise to a contradiction. Therefore  $d_{P_1}(v_{k-1})=0$. 
	Moreover, $N_H(v_{k-1})\subseteq V(P_1)\cup V(P_2)$, otherwise, $H-P_1$ has a path with at least $q+1$ vertices, which contradicts the choice of $P_2$.  Then $d_{H}(v_{k-1})=d_{P_2}(v_{k-1})$, and so
	$d_{H}(v_{k-1})=d_{P_2}(v_{k-1})=2$ by Lemma \ref{v2}. Clearly, $H-v_{k-1}$ is connected.
	Taking $w=v_{k-1}$,  (3) of the lemma occurs.   This completes the proof.
	\qed
	
	\vskip 10pt
	
	In what follows, we consider the existence of chorded cycles in a   $2$-connected graph.
	It is proved in \cite{BRK} that if a $2$-connected graph of order at least $4$
	contains no chorded cycles, then it is triangle-free. This gives rise to a sufficient condition for the existence of a chorded cycle in a $2$-connected graph.

	\begin{lemma}
		\label{triangle-free}
		Let $G$ be a $2$-connected graph of order at least $4$. If $G$ contains a triangle, then
		$G$ contains a chorded cycle.
	\end{lemma}
	
	In a $2$-connected  graph which is not a cycle, a longest cycle always has a good ear defined as follows.  
	
	\begin{definition}\label{good-ear}
		Let $G$ be a connected triangle-free graph, and let $C$ be a longest cycle in $G$, which has length $t$. Let $I$ be the vertex set of a path 
		on $C$, and let $\mathcal{E}_I$ be the set of ears of $C$ in $G$ each of which has two ends in $I$. Suppose that   $\mathcal{E}_I\ne\emptyset$. 
		Each member of $\mathcal{E}_I$ is called an $I$-ear of $C$ in $G$. An 
		$I$-ear $P$ of $C$  is said to be good if $P$ meets the following conditions in order: 
		\begin{itemize}
			\item[\rm (1)] the  ends of  $P$ are  as close as possible on  $C$,  
			\item[\rm (2)] the length of $P$ is  as large as possible.
		\end{itemize} 
	\end{definition}

	From Definition \ref{good-ear} it follows that there exists a set $I$ with $|I| \leqslant \frac{t}{2}+2$ such that $G$ has a good $I$-ear.

	\begin{lemma}\label{Basic}
		Let $G$ be a  $2$-connected  triangle-free graph, and let $C$ be a longest cycle in $G$, which has length $t$.
		Suppose that $P$ is a good $I$-ear of $C$ which is described as in Definition {\rm \ref{good-ear}}. Without loss of generality, write $C=u_1u_2\cdots u_tu_1$, $P=u_1v_1\cdots v_{\ell}u_k$ and  $I=\{u_i\,:\, 1\leqslant i\leqslant \frac{t}{2}+2\}$,
		where $k\geqslant 2$, $\ell \geqslant 0$, and if  $\ell=0$ then  $P=u_1u_k$. Then $G$ contains  chorded cycles provided that
		one of the following holds:
		\begin{itemize}
			\item[\rm (1)] $C$ has an ear of length $1$;
			\item[\rm(2)]  $d_{G}(u_{2}) \geqslant 3$ or $d_{G}(u_{k-1}) \geqslant 3$;
			\item[\rm (3)] $d_{G}(u_{i}) \geqslant 3$ and $d_{G}(u_{i+1}) \geqslant 3$, for some $1\leqslant i\leqslant k-1$;
			\item[\rm (4)] $d_G(w,z)\geqslant 4$ for each pair of  distinct vertices $w,z\in \{u_i,v_j\,:\, 2\leqslant i\leqslant k-1, 1\leqslant j\leqslant \ell\}$ with $wz\not\in E(G)$.
		\end{itemize}
	\end{lemma}
	\proof  
	First, if an ear of $C$ has length $1$ then it is itself a  chord of $C$, and $G$ contains  a 
	chorded cycle. In view of this, we suppose next that every ear of $C$ has length at least $2$. In particular, by the choices of $C$ and $P$,  we have
	$3 \leqslant k\leqslant \frac{t}{2}+2$ and $1\leqslant \ell\leqslant k-2$.

	{\it Case} 1. Suppose that (2) holds. Without loss of generality, we let $d_{G}(u_{2}) \geqslant 3$, and pick a 
	vertex $x$ of $G$ in $ N_G(u_2)\setminus \{u_1,u_3\}$.
	Since $G$ is $2$-connected, $G-u_2$ is connected. Pick a shortest path $Q[x,y]$ in $G-u_2$ that connects $x$ and the cycle $C[u_k,u_1]\cup P$. 
	We claim that $Q[x,y]$ has no vertices lying on the path $C[u_1,u_k]$. Suppose the contrary, and let $u_i$ be the first (from $x$) common vertex of  $Q[x,y]$ and $C[u_1,u_k]$. Then we get an $I$-ear $u_2x\cup Q[x,u_i]$  with ends closer on $C$ than those of $P$, which is not the case. Therefore,  $V(Q[x,y])\cap V(C[u_1,u_k])=\emptyset$. In addition, if $y$ lies on $P$, then we get a similar contradiction. This allows us to let $y=u_j$ for some $k+1 \leqslant j \leqslant t$. Now we have a cycle $u_1 \cup P\cup u_ku_{k-1}\cdots u_2x\cup Q[x,u_j] \cup C[u_j,u_1]$, which has a chord $u_1u_2$.

	{\it Case} 2. Suppose that (3) holds. In view of Case 1, we let $3\leqslant i, i+1\leqslant  k-2$, and so $k\geqslant 6$. Pick $x_i\in N_G(u_i)\setminus\{u_{i-1},u_{i+1}\}$, 
	$x_{i+1}\in N_G(u_{i+1})\setminus\{u_i,u_{i+2}\}$, a shortest path 
	$Q[x_i,y_i]$ in $G-u_{i}$ that connects $x_i$ and the cycle $C[u_k,u_1]\cup P$, 
	and  a shortest path $R[x_{i+1},y_{i+1}]$ in $G-u_{i+1}$ that connects
	$x_{i+1}$ and the cycle $C[u_k,u_1]\cup P$. 
	If either  $Q[x_i,y_i]$ or $R[x_{i+1},y_{i+1}]$ has a vertex lying on the cycle $C[u_1,u_k] \cup P$, 
	then an argument similar to that in Case 1 implies that  $C$ has an $I$-ear  
	with ends closer on $C$ than those of $P$, a contradiction. 
	Thus, we may put $y_i=u_{i'}$ and $y_{i+1}=u_{j'}$ with $k+1 \leqslant i',j' \leqslant t$. In addition, if $Q[x_i,u_{i'}]$ and $R[x_{i+1},u_{j'}]$ have a common internal vertex, then a similar contradiction arises. Let  
	$T=C[u_{j'},u_{i'}]$ if $j'\leqslant i'$, otherwise, let $T$ be the reverse sequence of $C[u_{i'},u_{j'}]$ if $i' \leqslant j'$. 
	Then we get a cycle $u_1 \cup P\cup u_ku_{k-1}\cdots u_{i+1}x_{i+1}\cup R[x_{i+1},u_{j'}]\cup T\cup Q[u_{i'},x_i] \cup x_iu_iu_{i-1}\cdots u_2u_1$, which has a chord $u_iu_{i+1}$.

	{\it Case} 3. 
	Finally, suppose that (4) holds. By (2), we let  $d_{G}(u_{2}) =2=d_{G}(u_{k-1})$.  
	Since $G$ is triangle-free, $u_2v_1\not\in E(G)$.
	Then $d_{G}(u_{2},v_{1}) \geqslant 4$. Since $u_{1} \in N_{G}(u_{2}) \cap N_{G}(v_{1})$,  we have $d_{G}(v_{1}) \geqslant 3$. 
	Recall  that $1\leqslant \ell\leqslant k-2$. 
	
	Suppose that $\ell=1$. Recalling that
	$d_{G}(v_{1}) \geqslant 3$, pick $x\in N_G(v_1)\setminus \{u_1,u_k\}$ and a shortest path $Q[x,y]$ in $G-v_1$ that connects $x$ and the cycle $C$. If $y$ lies on the path $C[u_1,u_k]$, then either
	$G$ has a triangle or $C$ has an $I$-ear with ends closer on $C$ than those of $P$, a contradiction.
	Thus, we may put $y=u_j$ for some $k+1 \leqslant j \leqslant t$. Then we have a cycle $u_1 \cup C[u_1,u_k]\cup u_kv_1x\cup Q[x,u_j]\cup C[u_j,u_1]$, which has a chord $u_1v_1$.

	Suppose that $k\geqslant 5$.
	Recalling  that $C$ has no ears of length $1$,  we have $u_iu_{i+2}\not\in E(G)$ for all $1\leqslant i\leqslant k-2$. Since $d_G(u_i,u_{i+2})\geqslant 4$, if $d_{G}(u_{j})= 2$ for some $1\leqslant j\leqslant k$, then either $d_{G}(u_{j+2})\geqslant 3$ or $d_{G}(u_{j-2})\geqslant 3$.  Now, since $d_{G}(u_{2})=2$, we have  $d_{G}(u_{4})\geqslant 3$, and $k\geqslant 6$ as $d_{G}(u_{k-1})= 2$. In addition, if $d_{G}(u_{3})=2$ then  $d_{G}(u_{5}) \geqslant 3$, and $k\geqslant 7$ as $d_{G}(u_{k-1})= 2$. These observations show that  there is $i$ with $3\leqslant i<i+1\leqslant k-3$ such that $d_{G}(u_{i})\geqslant 3$ and $d_{G}(u_{i+1})\geqslant 3$. Then (3) holds, and so 
	$G$ contains a chorded cycle. 
	By the argument above, we let $k=5$  and suppose that $d_G(u_2)=2$.  
	Since $u_2u_4\not\in E(G)$, we have $d_{G}(u_2,u_4)\geqslant 4$.
	It follows that $d_G(u_4)\geqslant 3$. Then (2)  of Lemma \ref{Basic} holds, and $G$ contains a chorded cycle.

	To complete the proof, we may let $k=4$ and $\ell=2$. 
	By (2), we get $d_G(u_2)=2=d_G(u_3)$, which implies $d_G(v_1) \geqslant 3$ and $d_G(v_2) \geqslant 3$.
	Since $d_G(v_1) \geqslant 3$, picking $x\in N_G(v_1)\setminus \{u_1,v_2\}$ and a shortest path $Q[x,y]$ in $G-v_1$ that connects $x$ and the cycle $C$, we have $y \in C[u_4,u_t]$ from a similar discussion above and let $y=u_j$ for some $4 \leqslant j \leqslant t$.
	If $5 \leqslant j \leqslant t$ then a cycle $v_1x\cup Q[x,u_j]\cup u_ju_{j+1} \cdots u_tu_1u_2u_3u_4v_2v_1$ with a chord $u_1v_1$ is obtained. Thus $j=4$. 
	Considering $d_G(v_2) \geqslant 3$ and choosing $x'\in N_G(v_2)\setminus \{v_1,u_4\}$ and a shortest path $R[x',y']$ in $G-v_2$ that connects $x'$ and the cycle $C$, we get $y' \in C[u_4,u_1]$ and let $y'=u_{j'}$ for some $j' \in \{1,5,\cdots,t\}$. Similar to the above discussion, we get a chorded cycle when $j' \in \{5,\cdots,t\}$. Now for $j'=1$, we get a cycle $v_1x\cup Q[x,u_4]\cup v_2x'\cup R[x',u_1]\cup u_1v_1$ with a chord $v_1v_2$. 
	This completes the proof.
	\qed

	\begin{corollary}\label{Basic-cor}
		Let $G$ be a  $2$-connected graph of order at least $4$ and  $\delta_{2}(G) \geqslant 4$. Then $G$ contains  chorded cycles. 
	\end{corollary}
	\proof 
	By Lemma \ref{triangle-free},   we may suppose that $G$ is triangle-free. 
	Let  $C=u_1u_2\cdots u_tu_1$ be a longest cycle in $G$, and so $t\geqslant 4$. Since $\delta_{2}(G) \geqslant 4$, we have $G\ne C$, and so $C$ has at least one ear in  $G$. Then the result follows from Lemma \ref{Basic}
	\qed
	
	\vskip 10pt

	A  {\em block} in a graph is a maximal subgraph without cut-vertices. 
	Recall that the blocks of a connected graph fit together in a tree-like
	structure. In particular, if a graph $G$ of order at least $3$ is connected but not $2$-connected,
	then $G$ has at least two blocks each of which contains a unique cut-vertex of $G$.
	For convenience, we call a block of a connected graph a  {\em leaf block} if it contains a unique cut-vertex of the graph. 
	The following result says
	that Corollary \ref{Basic-cor} holds for a connected graph that has no triangle blocks. 
	
	\begin{lemma}\label{3-neighbor}
		Let $G$ be a connected graph of order at least $4$ and $\delta_2(G)\geqslant 4$. 
		Then either $G$ contains a chorded cycle, or all leaf blocks of $G$ are triangles.  
	\end{lemma}
	\proof If $G$ is $2$-connected then the result is true by Corollary \ref{Basic-cor}.
	Suppose next that $G$ is not $2$-connected. Let $L_0$, $L_1$, $L_2,\ldots L_{m}$ be the leaf blocks of $G$, and let $x_i$ be the cut-vertex of $G$ in $L_i$, where $0\leqslant i\leqslant m$. We have $m\geqslant 1$. Clearly, if some $L_i$ is a complete graph of order at least $4$, then
	$G$ contains a chorded cycle. In addition,  by Corollary \ref{Basic-cor}, if $\delta_2(L_i)\geqslant 4$ for some $i$ then   $G$ contains a chorded cycle. 
	Thus we suppose further that for each $0\leqslant i\leqslant m$, neither $\delta_2(L_i)\geqslant 4$ nor 
	$L_i$ is a complete graph of order at least $4$.
	
	Clearly, $d_{G}(u)=d_{L_i}(u)$, $d_{G}(v)=d_{L_j}(v)$ and 
	$N_G(u,v)=N_{L_i}(u)\cup N_{L_j}(v)$ for $1\leqslant i,j\leqslant m$ with $u\in V(L_i)\setminus \{x_i\}$ and $v\in V(L_j)\setminus \{x_j\}$. In view of this, if
	every $L_i$ has at most three vertices then  every $L_i$ is a triangle, and if 
	some $L_i$ has order at least $4$ then $L_i$ is not a cycle.
	Thus we next suppose that one of the $L_i$, say $L_0$ without loss of generality, has order at least $4$.
	Put $B=L_0$. 
	
	Let $C$ be a longest cycle in $B$ of length say $t$.  Employing  Lemma \ref{Basic}, we next show that  $B$   
	contains a chorded cycle, and  the lemma follows. This is obvious when $C$ has an ear of length $1$.
	We suppose further that  $C$ has no ears of length $1$. In particular, 
	$t\geqslant 4$.
	Next, we continue the argument in two cases: $x_0\in V(C)$, and $x_0\not\in V(C)$.
	
	{\it Case} 1. Suppose that $x_0\in V(C)$. Let $x\in V(C)$ be at distance $\left[\frac{t}{2}\right]$ on $C$ from $x_0$, and put $N_C(x)=\{y_1,y_2\}$. Then $y_1y_2\not\in E(G)$ as  $C$ has no ears of length $1$.
	Since $C$ has length $t\geqslant 4$, we have $x_0\not\in \{y_1,y_2\}$, and hence  $d_{B}(y_1,y_2)=d_{G}(y_1,y_2)\geqslant 4$. Thus, without loss of generality, we let $d_B(y_1)\geqslant 3$.
	Again, since $C$ has no ears of length $1$, pick $w\in N_B(y_1)\setminus N_C(y_1)$.
	Considering the shortest path $Q[w,z]$ from $w$ to $C-y_1$, where $z \in V(C) \setminus \{y_1\}$, we may obtain a path $C[y_1,z]$ or $C[z,y_1]$  on $C$ of length $\left[\frac{t+1}{2}\right]$ with $x_0$ not an internal vertex of the path, and an ear of $C$ in $B$ with two ends lying on $C[y_1,z]$ or $C[z,y_1]$.  Let $I$ be the vertex set of this $\left[\frac{t+1}{2}\right]$-path. Then
	$\mathcal{E}_I\ne\emptyset$. Choose a good $I$-ear $P$ of $C$, see Definition \ref{good-ear}.
	It is easily shown that one of (2)-(4) of Lemma \ref{Basic} holds for the triple
	$(B,C,P)$. Then $B$ contains a chorded cycle, and the lemma follows.

	{\it Case} 2. Suppose that $x_0\not\in V(C)$. Without loss of generality, write $C=u_1u_2\cdots u_tu_1$. Take a   good $I$-ear, say $P=u_1v_1\cdots v_{\ell}u_k$, such that $x_0 \notin V(P)\setminus \{u_1,u_k\}$ as much as possible, where $I=\{u_i\,:\, 1\leqslant i\leqslant t-1\}$, $k \geqslant 3$ and $\ell \geqslant 1$. 
	If $x_0\not\in V(P) \setminus \{u_1,u_k\}$ then one of (2)-(4) of Lemma \ref{Basic} holds for the triple
	$(B,C,P)$, and the lemma follows. 
	Suppose next that  $x_0 \in V(P)\setminus \{u_1,u_k\}$. Of course,
	$x_0$ is an internal vertex of $P$, and hence $1\leqslant \ell\leqslant k-2$.
	
	Considering $d_{B}(u_{k-1},u_{k+1})=d_{G}(u_{k-1},u_{k+1})\geqslant 4$, we have either $d_B(u_{k-1})\geqslant 3$ or $d_B(u_{k+1})\geqslant 3$.
	Suppose first that $d_B(u_{k-1})\geqslant 3$. Pick $x\in N_G(u_{k-1})\setminus \{u_{k-2},u_k\}$ and a shortest path $Q[x,y]$ in $G-u_{k-1}$ that connects $x$ and the cycle $C$. Obviously $y$ does not lie on the path $C[u_1,u_k]$ and $Q[x,y] \cap P = \emptyset$.
	Thus, we may put $y=u_j$ for some $k+1 \leqslant j \leqslant t$. Then $u_{k-1}x \cup Q[x,u_j]$ is an ear of $C$ and does not contain $x_0$.
	Suppose next that $d_B(u_{k+1})\geqslant 3$. Pick $x'\in N_G(u_{k+1})\setminus \{u_{k},u_{k+2}\}$ and a shortest path $R[x',y']$ in $G-u_{k+1}$ that connects $x'$ and the cycle $C$. Obviously, $R[x',y'] \cap P = \emptyset$. Put $y'=u_{j'}$ for some $j' \notin \{k,k+1,k+2\}$. Then $u_{k+1}x' \cup R[x',u_{j'}]$ is an ear of $C$ and does not contain $x_0$.
	In the two cases mentioned above, choosing a good ear $P'$ not containing $x_0$ for $C$ and 
	repeating the argument in Case 1 for $(B, C, P')$, it follows that $B$ contains a chorded cycle, and the lemma follows.
	This completes the proof.
	\qed
	
	\vskip 10pt
	
	It is easy to deduce the following result from the proof of Lemma \ref{3-neighbor}.
	
	\begin{corollary}\label{3-neighbor-cor}
		Let $G$ be a $2$-connected graph of order at least $4$. Suppose that 
		$G$ has at most one vertex $x$ with the  property that $d_G(x,y)\leqslant 3$ for some
		$y\in V(G)\setminus (N_G(x)\cup\{x\})$.  
		Then  $G$ contains a chorded cycle.  
	\end{corollary}
	
	\vskip 20pt

	\section{Optimal systems of
		chorded cycles}

	For  a collection  $\mathcal{C}$ of  subgraphs in a graph $G$, we put $V(\mathcal{C})=\cup_{H\in \mathcal{C}}V(H)$ and $G-\mathcal{C}=\langle V(G)-V(\mathcal{C})\rangle$, where $G-\mathcal{C}$ is the null graph when $V(G)=V( \mathcal{C})$. Let $r$ be a positive integer. We call $\mathcal{C}$ a {\em minimal $r$-system} if 
	$|\mathcal{C}|=r$, $|V(\mathcal{C})|$ is as small as possible, and $\mathcal{C}$ consists of pairwise vertex-disjoint subgraphs.

	\begin{lemma}\label{C-mini} 
		Let $\mathcal{C}$ be a minimal  $r$-system  of  chorded cycles in a graph $G$, and $C\in \mathcal{C}$.
		Then $\langle V(C)\cup S\rangle$ contains no chorded cycles of length less than $|V(C)|$, where $S\subseteq V(G-\mathcal{C})$.
	\end{lemma}
	\proof Suppose the contrary that $\langle V(C)\cup S\rangle$ contains a chorded cycle $C'$ with
	$|V(C')|<|V(C)|$, where  $S\subseteq V(G-\mathcal{C})$.
	Then we have a collection $\mathcal{C}'$ of vertex-disjoint chorded cycles, which is obtained from $\mathcal{C}$ by replacing $C$ with $C'$. Clearly,  $r=|\mathcal{C}|= |\mathcal{C}'|$, but
	$|V(\mathcal{C}')|=|V(\mathcal{C})|-|V(C)|+|V(C')|<|V(\mathcal{C})|$, contrary to the hypothesis.
	\qed
	
	\begin{lemma}
		\label{degree-3}
		Let $\mathcal{C}$ be a minimal  $r$-system  of  chorded cycles in a graph $G$, and $C\in \mathcal{C}$. Suppose that $d_{C}(u) \geqslant 3$ for some $u\in V(G-\mathcal{C})$. Then 
		one of the following holds.   
		\begin{itemize}
			\item[\rm (1)]   $|V(C)|=4$, and $d_{C}(u) \in\{3,4\}$.
			\item[\rm (2)]  $|V(C)|=5$, $d_{C}(u)=3$, and two vertices in $N_C(u)$  are both at distance $2$ on $C$ from the third vertex in  $N_C(u)$.
			\item[\rm (3)] $|V(C)|=6$, $d_{C}(u)=3$, and  $\langle V(C)\cup \{v\}\rangle$ is triangle-free for any $v\in V(G-\mathcal{C})$.
		\end{itemize}
	\end{lemma}
	\proof 
	Write $C=u_1u_2\cdots u_tu_1$. Since $C$ is a chorded cycle, $t\geqslant 4$.  
	Choose  $v,w \in N_C(u)$ such that $v$ and $w$ are at 
	distance on $C$   as large as possible.
	Without loss of generality, assume that $v=u_1$, and $w=u_k$ with
	$2\leqslant k\leqslant \left[\frac{t+2}{2}\right]$. 
	If $k=2$ then, since $d_{C}(u) \geqslant 3$, it is easily deduced that  $C$ is a $3$-cycle, which is not the case.
	Therefore, $k\geqslant 3$.

	Suppose first that $N_C(u)$ contains some internal vertex of $C[u_1,u_k]$.
	By Lemma \ref{vertex-disjoint-paths}, the subgraph $C[u_1,u_k]\sqcup u$ contains a chorded cycle $C'$ of length no more than $k+1$. By Lemma \ref{C-mini}, $k+1\geqslant t$, and
	$t-1\leqslant k\leqslant \left[\frac{t+2}{2}\right]$, yielding $t=4$. It follows that $d_{C}(u) \in\{3,4\}$, as claimed in (1) of the lemma.
	
	Suppose now that $N_C(u)$ contains no internal vertices of $C[u_1,u_k]$.
	If $k\geqslant 4$ then the subgraph $C[u_k,u_1]\sqcup u$ contains a chorded cycle of length
	at most $t-1$, which contradicts Lemma \ref{C-mini}. We have $k=3$. This says that any two distinct vertices in $N_C(u)$ are at distance $1$ or $2$ on $C$. We deduce that $d_C(u)=3$, and 
	either  $|V(C)|=5$ and $N_C(u)=\{u_1,u_3,u_4\}$, or  $|V(C)|=6$ and $N_C(u)=\{u_1,u_3,u_5\}$. If $|V(C)|=5$ then we get (2)   of the lemma. 
	
	Suppose that  $|V(C)|=6$. If $C$ has an ear of length $1$ say $u_1u_3$ or $u_2u_4$, for example, then we have a $4$-cycle $u_1uu_3u_2u_1$ with a chord $u_1u_3$ or a $5$-cycle 
	$u_2u_3uu_5u_4u_2$ with a chord $u_3u_4$, which contradicts Lemma \ref{C-mini}.
	This says that $\langle V(C)\rangle$ is triangle-free, in particular, each chord
	of $C$ joins two antipodal vertices on $C$. Without loss of generality, let $u_1u_4$ be a chord of $C$.
	Suppose that $\langle V(C)\cup \{v\}\rangle$ contains a triangle for some  $v\in V(G-\mathcal{C})$.
	Without loss of generality, let $\{u_1,u_2\}\subseteq N_C(v)$ or  $\{u_2,u_3\}\subseteq N_C(v)$ or $\{u_1,u_4\}\subseteq N_C(v)$.
	Then $u_1vu_2u_3u_4u_1$ is a $5$-cycle with a chord $u_1u_2$ or $u_1u_2vu_3u_4u_1$ is a $5$-cycle with a chord $u_2u_3$ or $u_1vu_4u_3u_2u_1$ is a $5$-cycle with a chord $u_1u_4$, which contradicts Lemma \ref{C-mini}. Thus  (3)  of the lemma follows. 
	\qed

	\begin{lemma}\label{6-cycle}
		Let $\mathcal{C}$ be a minimal  $r$-system  of  chorded cycles in a graph $G$, and $C\in \mathcal{C}$.
		Suppose that $d_{C}(u,v) \geqslant 5$ for some  $u,v \in V(G-\mathcal{C})$.
		Then $|V(C)|=6$, and there exist $u'\in N_C(u)$ and $v'\in N_C(v)$ such that
		$C_v=(C-u')\sqcup v$ and $C_u=(C-v')\sqcup u$ are chorded $6$-cycles.
	\end{lemma}
	\proof
	
	Since  $d_{C}(u,v)  \geqslant 5$, we have $|V(C)| \geqslant 5$ and, without loss of generality, let
	$d_C(u)\geqslant 3$. 
	By Lemma \ref{degree-3},  $d_C(u)=3$ and  $|V(C)|=5$ or $6$.
	Of course, $d_C(v)\geqslant 2$.  
	Write $C=u_1u_2\cdots u_tu_1$, and suppose that $u_1\in N_C(u)$.
	
	Suppose that $t=5$. Then $N_C(u)=\{u_1,u_3,u_4\}$, and $\{u_2,u_5\}\subseteq N_C(v)$. It is easy to check that there exists a chorded cycle of length $4$ with vertices in $V(C)\cup\{u,v\}$, which contradicts Lemma \ref{C-mini}. 
	
	Then $t=6$  and, by Lemma \ref{degree-3}, $N_C(u)=\{u_1,u_3,u_5\}$, and each chord of $C$ joins two antipodal vertices on $C$. Note that $|N_C(v)\cap \{u_2,u_4,u_6\}|\geqslant 2$ and let $\{u_2,u_4\}\subseteq N_C(v)\subseteq \{u_2,u_4,u_6\}$. 
	First, observe that $uu_1u_6u_5u_4u_3u$ is a $6$-cycle with a chord $uu_5$, and $C_u=(C-v')\sqcup u$ is a chorded $6$-cycle by taking $v'=u_2$.
	Now if $u_1u_4$ or $u_2u_5$ is a chord of $C$ then $vu_2u_1u_6u_5u_4v$ is a $6$-cycle with a chord $u_1u_4$ or $u_2u_5$, and the lemma follows by taking $u'=u_3$.
	Thus, the only remaining possibility is that $u_3u_6$ is a chord of $C$. In this case,  $vu_2u_3u_6u_5u_4v$ is a  $6$-cycle with  chord $u_3u_4$. Taking $u'=u_1$ gives the desired conclusion. 
	\qed

	\begin{lemma}\label{5-path}
		Let $\mathcal{C}$ be a minimal  $r$-system  of  chorded cycles in a graph $G$, and $C\in \mathcal{C}$.
		Suppose that $H =G-\mathcal{C}$ contains a path $P=x_{1}x_{2} x_{3}x_{4}\cdots x_{\ell}$, where $ {\ell}\geqslant 4$.  Then 
		\begin{itemize}
			\item [\rm (1)] 
			$d_{C}(x_{1},x_{3})+d_{C}(x_{1},x_{4})+d_{C}(x_{2},x_{4})\leqslant 12$,
			\item [\rm (2)] 
			$d_{C}(x_{1},x_{3})+d_{C}(x_{2},x_{4})+d_{C}(x_{3},x_{5})\leqslant 12$,
			\item [\rm (3)] 
			$d_{C}(x_{1},x_{3})+d_{C}(x_{1},x_{4})+d_{C}(x_{2},x_{5})\leqslant 12$ if $x_2x_4\in E(G)$. 
			
		\end{itemize} 
		
	\end{lemma}
	
	\proof Our strategy is to investigate the subgraph $C\sqcup\langle  x_1,x_2,\ldots, x_k\rangle$, where $k\in\{4,5\}$. Suppose that the sum of three terms $d_C(x_i,x_j)$ in  (1), (2) or (3) is greater than $12$.
	Then at least one $d_C(x_i,x_j)$ of the three summands is not less than $5$.
	By Lemma \ref{6-cycle}, $|V(C)|=6$.
	By    Lemma \ref{degree-3},  
	for every $x_i$, the subgraph $C\sqcup x_i$ is triangle-free, and $d_{C}(x_{i})\leqslant 3$.
	In particular, we can assert that the following conclusions are valid:
	
	\begin{itemize}
		\item[(i)] if $d_C(x_i,x_j)\geqslant 5$ for distinct $i,j$, then each of $x_i$ and $x_j$ has at least two distinct neighbors on $C$, and, for each $x \in \{x_i,x_j\}$, any two neighbors of $x$ on $C$ are at distance $2$ on $C$;
		\item[(ii)] if $d_C(x_i) \geqslant 2$ for some $i$ and all neighbors of $x_i$ are at distance $2$ on $C$ from every other, then $x_i$ and $x_{i\pm 1}$ have no common neighbors  on $C$;
		\item[(iii)] if  $d_C(x_i,x_{i+2})\geqslant 5$ for some $i$, then
		$d_C( x_{i+1})\leqslant 1$ and $N_C(x_{i+1})\subseteq V(C)\setminus N_C(x_i,x_{i+2})$.
	\end{itemize}
	Based on these observations, we shall deduce the contradiction.
	Note that the positions of the two vertex pairs $(x_1,x_3)$ and $(x_2,x_4)$ on $P[x_1,x_4]$ are symmetrical,
	the positions of the two vertex pairs $(x_1,x_3)$ and $(x_3,x_5)$ on $P[x_1,x_5]$ are symmetrical,
	and the positions of the two vertex pairs $(x_1,x_4)$ and $(x_2,x_5)$ on $P[x_1,x_5]$ are symmetrical.
	We need only deal with the following three cases:  $d_{C}(x_{1},x_{3})\geqslant 5$ or  $d_{C}(x_{1},x_{4})\geqslant 5$ for (1),  $d_{C}(x_{1},x_{3})\geqslant 5$ or  $d_{C}(x_{2},x_{4})\geqslant 5$ for (2), and  $d_{C}(x_{1},x_{3})\geqslant 5$ or  $d_{C}(x_{1},x_{4})\geqslant 5$ for (3).
	We write $C=u_{1}u_{2}u_3u_4u_5u_{6}u_1$.

	{\it Case  }1. Suppose $d_{C}(x_{1},x_{3})\geqslant 5$ or  $d_{C}(x_{1},x_{4})\geqslant 5$ for (1).
	Without loss of generality, let $\{u_1,u_2,u_3,u_4,u_5\}\subseteq N_C(x_1,x_3)$ or  $\{u_1,u_2,u_3,u_4,u_5\}\subseteq N_C(x_1,x_4)$, respectively.
	
	Suppose first that   $d_{C}(x_{1},x_{3})\geqslant 5$.  Then, by (i) and (ii), $N_C(x_1,x_3)\cap N_C(x_2)=\emptyset$   and $N_C(x_3)\cap N_C(x_4)=\emptyset$.
	Suppose that  $d_{C}(x_{3})=3$. Then, by the assertion (i), we may let
	$N_C(x_3)=\{u_1,u_3,u_5\}$, and  so 
	$\{u_2,u_4\}\subseteq N_C(x_1)\subseteq\{u_2,u_4,u_6\}$, $N_C(x_2)\subseteq\{u_6\}$
	and $N_C(x_4)\subseteq\{u_2,u_4,u_6\}$.
	Then $d_C(x_1,x_3)\leqslant 6$, $d_C(x_1,x_4)\leqslant 3$ and $d_C(x_2,x_4)\leqslant 3$,
	yielding $13\leqslant d_{C}(x_{1},x_{3})+d_{C}(x_{1},x_{4})+d_{C}(x_{2},x_{4})
	\leqslant 6+3+3=12$, a contradiction. 
	Thus  $d_{C}(x_{1})=3$. Also by  the assertion  (i),  we may let
	$N_C(x_1)=\{u_1,u_3,u_5\}$. Then  $\{u_2,u_4\}\subseteq N_C(x_3)\subseteq\{u_2,u_4,u_6\}$, 
	$N_C(x_2)\subseteq\{u_6\}$ and $N_C(x_4)\subseteq\{u_1,u_3,u_5,u_6\}$.
	If   $N_C(x_3)=\{u_2,u_4,u_6\}$ then $N_C(x_2,x_4)\subseteq N_C(x_1)=\{u_1,u_3,u_5\}$,   and so  $13\leqslant d_{C}(x_{1},x_{3})+d_{C}(x_{1},x_{4})+d_{C}(x_{2},x_{4})
	\leqslant 6+3+3=12$, a contradiction. 
	Therefore,  $N_C(x_3)=\{u_2,u_4\}$.  We have $d_C(x_1,x_3)=5$, $d_C(x_1,x_4)\leqslant 4$ and $d_C(x_2,x_4)\leqslant 4$.
	Noting that $d_{C}(x_{1},x_{4})+d_{C}(x_{2},x_{4})\geqslant 13-d_{C}(x_{1},x_{3})=8$, it follows that
	$d_{C}(x_{1},x_{4})=d_{C}(x_{2},x_{4})=4$.
	By $d_{C}(x_{2},x_{4})=4$, we have   $N_C(x_2)=\{u_6\}$ and $N_C(x_4)= \{u_1,u_3,u_5\}$. 
	This yields that $d_{C}(x_{1},x_{4})=3$, a contradiction.

	Now let $d_{C}(x_{1},x_{4})  \geqslant  5$. Without loss of generality, by  the assertion (i), let $N_C(x_1)=\{u_1,u_3,u_5\}$  and  $N_C(x_4)=\{u_2,u_4,u_6\}$ or $\{u_2,u_4\}$. Using  the assertion (ii),  $N_C(x_1)\cap N_C(x_2)=\emptyset=N_C(x_3)\cap N_C(x_4)$. In particular, $ N_C(x_2)\subseteq\{u_2,u_4,u_6\}$ and  $ N_C(x_3)\subseteq\{u_1,u_3,u_5,u_6\}$. Then $d_C(x_1,x_3)\leqslant 4$, $d_C(x_1,x_4)\leqslant 6$ and 
	$d_C(x_2,x_4)\leqslant 3$. It follows from the hypothesis that
	$d_C(x_1,x_3)= 4$, $d_C(x_1,x_4)= 6$ and 
	$d_C(x_2,x_4)= 3$. By $d_C(x_1,x_4)= 6$, we have $N_C(x_4)=\{u_2,u_4,u_6\}$, and so
	$N_C(x_3)\subseteq\{u_1,u_3,u_5\}$ as $N_C(x_3)\cap N_C(x_4)=\emptyset$.
	Then $d_C(x_1,x_3)= 3$, a contradiction.

	{\it Case  }2. Suppose that $d_{C}(x_{1},x_{3})\geqslant 5$ or  $d_{C}(x_{2},x_{4})\geqslant 5$ for (2).
	
	{\it Subcase  }2.1.
	Suppose first that $d_{C}(x_{1},x_{3})  \geqslant 5$.  Then
	$N_C(x_1,x_3)\cap N_C(x_2)=\emptyset$ and $N_C(x_3)\cap N_C(x_4)=\emptyset$ by (i) and (ii).
	Without loss of generality, let $\{u_1,u_2,u_3,u_4,u_5\}\subseteq N_C(x_1,x_3)$.
	Then $N_C(x_2)\subseteq\{u_6\}$, see the assertion (iii).

	Suppose that  $d_{C}(x_{1})=3$, and let
	$N_C(x_1)=\{u_1,u_3,u_5\}$ without loss of generality, and  so 
	$\{u_2,u_4\}\subseteq N_C(x_3)\subseteq\{u_2,u_4,u_6\}$, and $N_C(x_4)\subseteq\{u_1,u_3,u_5,u_6\}$.
	In particular, $d_{C}(x_{1},x_{3})+d_{C}(x_{2},x_{4})\leqslant 9$.
	If $d_C(x_4)=0$  then, noting that $N_C(x_1,x_3)\cap N_C(x_2)=\emptyset$ and  $N_C(x_2)\subseteq\{u_6\}$, we have $13\leqslant  d_{C}(x_{1},x_{3})+d_{C}(x_{2},x_{4})+d_{C}(x_{3},x_{5})= (d_{C}(x_{1},x_{3})+d_{C}(x_{2}))+d_{C}(x_{3},x_{5}) \leqslant 6+6=12$, a contradiction. 
	If $d_C(x_4)=1$ then  $d_{C}(x_{1},x_{3})+d_{C}(x_{2},x_{4})\leqslant 7$, and so $13\leqslant  d_{C}(x_{1},x_{3})+d_{C}(x_{2},x_{4})+d_{C}(x_{3},x_{5})=7+d_{C}(x_{3},x_{5})$, yielding $d_{C}(x_{3},x_{5})\geqslant 6$. Thus $N_C(x_3)=\{u_2,u_4,u_6\}$ and $N_C(x_5)=\{u_1,u_3,u_5\}$, which contradicts $N_C(x_3,x_5)\cap N_C(x_4)=\emptyset$ by (iii).
	If $N_C(x_4)= \{u_1,u_3,u_5\}$ then $N_C(x_4)\cap N_C(x_5)=\emptyset$ by (ii), and $N_C(x_5)\subseteq \{u_2,u_4,u_6\}$ and $d_C(x_3,x_5) \leqslant 3$. Thus $13\leqslant  d_{C}(x_{1},x_{3})+d_{C}(x_{2},x_{4})+d_{C}(x_{3},x_{5}) \leqslant 9+3=12$, a contradiction.  
	Thus we have $d_C(x_4) =2$. This implies that either $N_C(x_4) = \{u_3,u_6\}$ or $N_C(x_4)\subseteq \{u_1,u_3,u_5\}$. For the former, we get $N_C(x_3)=\{u_2,u_4\}$, $d_C(x_2,x_4) = 2$ and $d_C(x_3,x_5) \leqslant 5$, and so $13\leqslant  d_{C}(x_{1},x_{3})+d_{C}(x_{2},x_{4})+d_{C}(x_{3},x_{5}) \leqslant 5+2+5=12$, a contradiction.  
	If $d_C(x_4)=2$ and $N_C(x_4)\subseteq \{u_1,u_3,u_5\}$, without loss of generality let $N_C(x_4)\subseteq \{u_1,u_3\}$, we have $N_C(x_4)\cap N_C(x_5)=\emptyset$ by (ii) and either $N_C(x_5)\subseteq \{u_2,u_4,u_6\}$ or $N_C(x_5)\subseteq \{u_2,u_5\}$. 
	For $N_C(x_5)\subseteq \{u_2,u_4,u_6\}$, we have $d_C(x_3,x_5) \leqslant 3$ and so $13\leqslant  d_{C}(x_{1},x_{3})+d_{C}(x_{2},x_{4})+d_{C}(x_{3},x_{5}) \leqslant 9+3=12$, a contradiction. 
	If $N_C(x_5) \subseteq \{u_2,u_5\}$ then $d_C(x_3,x_5) \leqslant 4$ and so $13\leqslant  d_{C}(x_{1},x_{3})+d_{C}(x_{2},x_{4})+d_{C}(x_{3},x_{5}) \leqslant 9+d_{C}(x_{3},x_{5}) \leqslant 13$. This implies that $d_C(x_3,x_5) =4$, $N_C(x_3)= \{u_2,u_4,u_6\}$, $\{u_5\} \subseteq N_C(x_5)\subseteq \{u_2,u_5\}$ and $d_C(x_2)=0$, and so $13\leqslant  d_{C}(x_{1},x_{3})+d_{C}(x_{2},x_{4})+d_{C}(x_{3},x_{5}) \leqslant 6+2+4=12$, a contradiction.

	Now let  $d_{C}(x_{3})=3$. Then, by the assertion (i), we may let
	$N_C(x_3)=\{u_1,u_3,u_5\}$, and  so $\{u_2,u_4\}\subseteq N_C(x_1)\subseteq\{u_2,u_4,u_6\}$, and $N_C(x_4)\subseteq\{u_2,u_4,u_6\}$. In particular, $d_{C}(x_{2},x_{4})\leqslant 3$. 
	If $d_C(x_4)=0$ then $d_{C}(x_{1},x_{3})+d_{C}(x_{2},x_{4})\leqslant 6$, and so $d_{C}(x_{3},x_{5})\geqslant 7$, a contradiction.
	Thus  $d_C(x_4) \geqslant 1$. 
	If $u_i \in N_C(x_4)$ then $u_i \notin N_C(x_5)$, otherwise $x_3x_4x_5u_iu_{i-1}x_3$ is a $5$-cycle with a chord $x_4u_i$, where $i \in \{2,4,6\}$, contrary to Lemma \ref{C-mini}. This implies that $d_C(x_3,x_5)+d_C(x_4) \leqslant 6$, and so $13\leqslant  d_{C}(x_{1},x_{3})+d_{C}(x_{2},x_{4})+d_{C}(x_{3},x_{5})=(d_{C}(x_{1},x_{3})+d_{C}(x_{2}))+(d_C(x_{4})+d_{C}(x_{3},x_{5})) \leqslant 6+6=12$, a contradiction.

	{\it Subcase  }2.2.
	Suppose  that $d_{C}(x_{2},x_{4})  \geqslant 5$.  
	Since the positions $x_2$ and $x_4$ on $P[x_1,x_5]$ are symmetric, we may let
	$d_{C}(x_{2})  \geqslant d_{C}(x_{4}) \geqslant 2$. 
	By the assertion (i), without loss of generality, we let $N_C(x_2)=\{u_1,u_3,u_5\}$, and
	so $\{u_2,u_4\}\subseteq N_C(x_4)\subseteq \{u_2,u_4,u_6\}$.
	Then  $N_C(x_1)\subseteq \{u_2,u_4,u_6\}$, and
	$N_C(x_3)\subseteq \{u_6\}$ with $d_C(x_3)+d_C(x_4)\leqslant 3$. In particular,
	$d_{C}(x_{1},x_{3})  \leqslant 3$.
	Then $13\leqslant  d_{C}(x_{1},x_{3})+d_{C}(x_{2},x_{4})+d_{C}(x_{3},x_{5})\leqslant 
	3+(d_{C}(x_{2})+d_C(x_{4})+d_{C}(x_{3}))+d_C(x_{5})\leqslant 9+d_C(x_{5})$,
	yielding $d_C(x_{5})\geqslant 4$,  a contradiction.

	{\it Case  }3. Suppose that $x_2x_4\in E(G)$, and
	$d_{C}(x_{1},x_{3})\geqslant 5$ or  $d_{C}(x_{1},x_{4})\geqslant 5$ for (3).
	Without loss of generality, let $\{u_1,u_2,u_3,u_4,u_5\}\subseteq N_C(x_1,x_3)$ or  $\{u_1,u_2,u_3,u_4,u_5\}\subseteq N_C(x_1,x_4)$, respectively.
	In addition, for any distinct $i,j \in \{2,3,4\}$, if $N_C(x_i)\cap N_C(x_j)\ne\emptyset$ then $G$ contains a $4$-cycle with a chord $x_ix_j$, contrary to Lemma \ref{C-mini}. Thus $N_C(x_i)\cap N_C(x_j)=\emptyset$.
	
	{\it Subcase  }3.1.
	Suppose   that $d_{C}(x_{1},x_{3})  \geqslant 5$.  Then
	$N_C(x_2)\subseteq \{u_6\}$. 
	
	Suppose first that $d_{C}(x_{1})=3$. Then $N_C(x_1)=\{u_1, u_3, u_5\}$, and
	$\{u_2, u_4\}\subseteq N_C(x_3)\subseteq \{u_2, u_4, u_6\}$ with $d_C(x_2)+d_C(x_3)\leqslant 3$.
	Since $N_C(x_3)\cap N_C(x_4)=\emptyset$,
	either $N_C(x_4)\subseteq \{u_3, u_6\}$ or $N_C(x_4)\subseteq \{u_1, u_3, u_5\}$.
	If $d_C(x_4)=0$ or $N_C(x_4)\subseteq \{u_1, u_3, u_5\}$, then $d_{C}(x_{1},x_{4})=3$, and so
	$13\leqslant  d_{C}(x_{1},x_{3})+d_{C}(x_{1},x_{4})+d_{C}(x_{2},x_{5})\leqslant 
	d_C(x_1)+d_C(x_3)+3+d_C(x_2)+d_C(x_5)\leqslant 3+3+3+d_C(x_5)$,  yielding $d_C(x_{5})\geqslant 4$,  a contradiction. 
	If $u_3 \in N_C(x_4)$ then $x_1x_2x_3x_4u_3x_1$ is a $5$-cycle with a chord $x_2x_4$, contrary to Lemma \ref{C-mini}.
	Thus we have  $N_C(x_4)=\{u_6\}$ then $N_C(x_3)=\{u_2, u_4\}$. 
	Recalling that  $N_C(x_2)\subseteq \{u_6\}$ and  $N_C(x_2)\cap N_C(x_4)=\emptyset$, we have 
	$d_C(x_2)=0$. Then  $13\leqslant  d_{C}(x_{1},x_{3})+d_{C}(x_{1},x_{4})+d_{C}(x_{2},x_{5}) \leqslant  
	5+4+3=12$, a contradiction.

	Now let  $d_{C}(x_{3})=3$.   Then $N_C(x_3)=\{u_1, u_3, u_5\}$,  
	$\{u_2, u_4\}\subseteq N_C(x_1)=\{u_2, u_4, u_6\}$ with $d_C(x_1)+d_C(x_2)\leqslant 3$,
	and $N_C(x_4)\subseteq  \{u_2, u_4, u_6\}$. In particular, $d_C(x_1,x_4)\leqslant 3$. Then
	$13\leqslant  d_{C}(x_{1},x_{3})+d_{C}(x_{1},x_{4})+d_{C}(x_{2},x_{5})
	\leqslant d_C(x_1)+d_C(x_3)+d_{C}(x_{1},x_{4})+d_{C}(x_{2})+d_{C}(x_{5})
	\leqslant  d_C(x_1)+d_{C}(x_{2})+ 3+3+d_{C}(x_{5})\leqslant 9+d_{C}(x_{5})$, yielding 
	$d_{C}(x_{5})\geqslant 4$, a contradiction.

	{\it Subcase  }3.2.
	Suppose   that $d_{C}(x_{1},x_{4})  \geqslant 5$.  
	Recall that $\{u_1,u_2,u_3,u_4,u_5\}\subseteq N_C(x_1,x_4)$.
	
	Suppose first that $d_{C}(x_{1})  =3$.  
	Then $N_C(x_1)=\{u_1, u_3, u_5\}$, and
	$\{u_2, u_4\}\subseteq N_C(x_4)\subseteq \{u_2, u_4, u_6\}$.
	Noting that 
	$N_C(x_1)\cap N_C(x_2)=N_C(x_2)\cap N_C(x_4)=N_C(x_3)\cap N_C(x_4)=N_C(x_4)\cap N_C(x_5)=\emptyset$,
	it follows that  $N_C(x_2)\subseteq \{u_6\}$ with $d_C(x_2)+d_C(x_4)\leqslant 3$, and
	either $N_C(x_3)\subseteq \{u_1, u_3, u_5\}$ or
	$d_C(x_3)+d_C(x_4)\leqslant 4$ with  $N_C(x_3)\subseteq \{u_3,u_6\}$.
	If $N_C(x_3)\subseteq \{u_1, u_3, u_5\}$, then
	$13\leqslant  d_{C}(x_{1},x_{3})+d_{C}(x_{1},x_{4})+d_{C}(x_{2},x_{5})
	\leqslant  d_{C}(x_{1},x_{3})+  d_{C}(x_{1})+  d_{C}(x_{4})+d_{C}(x_{2})+d_{C}(x_{5})
	\leqslant  9+d_{C}(x_{5})$, yielding 
	$d_{C}(x_{5})\geqslant 4$, a contradiction.
	If $u_3 \in N_C(x_3)$ then $x_3x_2x_4u_4u_3x_3$ is a $5$-cycle with a chord $x_3x_4$, contrary to Lemma \ref{C-mini}.
	This forces  $d_C(x_3)+d_C(x_4)\leqslant 3$ and $N_C(x_3)\subseteq \{u_6\}$.
	We have $d_{C}(x_{1},x_{3})+d_{C}(x_{1},x_{4})\leqslant 9$, and so $d_{C}(x_2,x_{5})\geqslant 4$.
	Recalling that $N_C(x_2)\subseteq \{u_6\}$ and $d_C(x_2)+d_C(x_4)\leqslant 3$, we have
	$N_C(x_2)=\{u_6\}$,   $N_C(x_5)=\{u_1, u_3, u_5\}$, and $N_C(x_4)=\{u_2, u_4\}$.
	Then $13\leqslant  d_{C}(x_{1},x_{3})+d_{C}(x_{1},x_{4})+d_{C}(x_{2},x_{5})= d_{C}(x_{1},x_{3})+5+4$,
	yielding $d_{C}(x_{1},x_{3})\geqslant 4$. We have $d_C(x_3) \geqslant 1$, and so  $N_C(x_3)= \{u_6\}$.
	Then $u_6x_2x_4x_3u_6$ is a $4$-cycle  with a chord $x_2x_3$,
	contrary to Lemma \ref{C-mini}. 
	
	Now let $d_{C}(x_{4}) =3$.  Then $N_C(x_4)=\{u_1, u_3, u_5\}$,  
	$\{u_2, u_4\}\subseteq N_C(x_1)\subseteq \{u_2, u_4, u_6\}$,  $N_C(x_3,x_5)\subseteq\{u_2, u_4, u_6\}$,
	and either $N_C(x_2)\subseteq\{u_6\}$ with $d_C(x_1)+d_C(x_2)\leqslant 3$ or  $N_C(x_2)\subseteq\{u_1, u_3, u_5\}$. For $N_C(x_2)\subseteq\{u_6\}$, we have $13\leqslant  d_{C}(x_{1},x_{3})+d_{C}(x_{1},x_{4})+d_{C}(x_{2},x_{5})
	\leqslant  3+6+3=12$, a contradiction. Thus $N_C(x_2)\subseteq\{u_1, u_3, u_5\}$, and so
	$d_C(x_2)=0$ as $N_C(x_4)\cap N_C(x_2)=\emptyset$.
	Then  $13\leqslant  d_{C}(x_{1},x_{3})+d_{C}(x_{1},x_{4})+d_{C}(x_{2},x_{5})
	\leqslant  3+6+3=12$, a contradiction.
	\qed

	\vskip 10pt

	For a collection  $\mathcal{C}$ of subgraphs in a graph $G$, denote by $r(G,\mathcal{C})$ the order   of a component in $G-\mathcal{C}$ with maximal order. Note that if $V(G)=V( \mathcal{C})$ then we put  $r(G,\mathcal{C})=0$.

	\begin{definition}\label{D3}
		A minimal  $r$-system $\mathcal{C}$ of  chorded cycles in a graph $G$ is called an optimal $r$-system of  chorded cycles if $r(G,\mathcal{C})$ is as large as possible.
	\end{definition}
	
	Let $G$ be a $2$-connected graph of order at least $4s$ and $\delta_{2}(G) \geqslant 4s$, where $s \geqslant 2$. Pick $S\subset V(G)$ with $|S|\leqslant 3$, and consider the graph $G-S$.  
	Then, for $u,v\in V(G-S)$ with $uv\not\in E(G-S)$, we have $uv\not\in E(G)$, and
	$d_{G-S}(u,v) \geqslant 4s-3=4(s-1)+1$. By Theorem \ref{T2}, $G-S$ and hence $G$ contains
	$s-1$ vertex-disjoint chorded cycles. Thus we may choose in $G$ an optimal $(s-1)$-system of  chorded cycles.

	\begin{lemma}\label{4-path-cycle}
		Let $G$ be a  $2$-connected graph of order at least $4s$ and $\delta_{2}(G) \geqslant 4s$, where $s \geqslant 2$. Let $\mathcal{C}$ be an optimal $(s-1)$-system of chorded cycles in $G$, and let $H$ be a component of order $4$ in $G-\mathcal{C}$. Then $H$ has no Hamiltonian paths.
	\end{lemma}
	\proof 
	Suppose that $H$ is a path. Write $H=uu'vv'$.
	Then $d_H(u,v)+d_H(u',v')+d_H(u,v')\leqslant 6$, and so
	$\sum_{C\in \mathcal{C}} (d_C(u,v)+d_C(u',v')+d_C(u,v'))\geqslant 12s-6$.
	Thus $d_C(u,v)+d_C(u',v')+d_C(u,v')\geqslant 13$ for some $C\in \mathcal{C}$, which contradicts
	(1) of Lemma \ref{5-path}. 
	
	Suppose  that $H$ is a $4$-cycle, and write $H=uu'vv'u$.
	Then $d_H(u,v)+d_H(u',v')\leqslant 4$, and  $\sum_{C\in \mathcal{C}}(d_C(u,v)+d_C(u',v'))\geqslant 8s-4$. Pick $C\in \mathcal{C}$ with  $d_C(u,v)+d_C(u',v')\geqslant 9$  and, without loss of generality,  let $d_C(u,v)\geqslant 5$. Then $|V(C)|=6$, and $d_C(u',v')\geqslant 3$.
	It follows that either $u'$ or $v'$ shares a neighbor on $C$ with one of $u$ and $v$.
	This shall give rise to a chorded cycle of length $5$, which contradicts Lemma \ref{C-mini}.
	
	Suppose that $H$ is a triangle plus a hanging edge, which has vertex set $\{u,u',v,v'\}$ and edge set $\{uv, uu',u'v,vv'\}$.
	Then $d_H(u,v')+d_H(u',v')\leqslant 4$, and  $\sum_{C\in \mathcal{C}}(d_C(u,v')+d_C(u',v'))\geqslant 8s-4$. Pick $C\in \mathcal{C}$ with  $d_C(u,v')+d_C(u',v')\geqslant 9$  and, without loss of generality,  let $d_C(u,v')\geqslant 5$. Then $|V(C)|=6$, write $C=w_1w_2w_3w_4w_5w_6w_1$.
	
	If $d_C(u)=3$, let $N_C(u)=\{w_1,w_3,w_5\}$, then without loss of generality let $w_2 \in N_C(v')$.
	Note that $u'w_2 \notin E(G)$ and $vw_2 \notin E(G)$, otherwise $H \sqcup w_2$ has a chorded $5$-cycle, which contradicts Lemma \ref{C-mini}.
	Replace $C$ with the new chorded $6$-cycle $C'=uw_1w_6w_5w_4w_3u$ with a chord $uw_5$, and 
	$G-(\mathcal{C}\setminus C\cup C')$ is a $4$-path $u'vv'w_2$. 
	
	Thus $d_C(u)=2$ and let $N_C(u)=\{w_1,w_3\}$. Obviously, $N_C(v')=\{w_2,w_4,w_6\}$ and $u'w_2, vw_2,u'w_4, vw_4 \notin E(G)$. If $w_1w_4 \in E(G)$ or $w_3w_6 \in E(G)$ then let $C'=uw_1w_6w_5w_4w_3u$ with a chord $w_1w_4$ or $w_3w_6$.  Thus $w_2w_5 \in E(G)$ and then pick $C'=uw_3w_2w_5w_6w_1u$ with a chord $w_2w_5$.  Then replace $C$ with the new chorded $6$-cycle $C'$, and clearly $G-(\mathcal{C}\setminus C\cup C')$ is a $4$-path. 
	By applying an analogous argument from the first paragraph to the two cases mentioned above, we derive a contradiction.
	Then the lemma follows.
	\qed
	
	\begin{lemma}\label{triangle-leaf-block}
		Let $G$ be a $2$-connected graph of order at least $4s$ and $\delta_{2}(G) \geqslant 4s$, where $s \geqslant 2$. Let $\mathcal{C}$ be an optimal $(s-1)$-system of chorded cycles in $G$, and let $H$ be a component of maximal order in $G-\mathcal{C}$.
		Suppose that $H$  contains two triangle leaf blocks.
		Then $G$ contains a collection of $s$ vertex-disjoint chorded cycles.
	\end{lemma}
	\proof 
	Let $B_x$ and $B_y$ be two triangle leaf blocks in $H$ containing the cut-vertices $x$ and $y$, respectively. Write $V(B_x)=\{x,x_1,x_2\}$ and $V(B_y)=\{y,y_1,y_2\}$.
	By Lemma \ref{degree-3}, $d_{C}(x_i)\leqslant 4$ and $d_{C}(y_j)\leqslant 4$, where $i,j\in \{1,2\}$ and $C \in \mathcal{C}$.
	Considering the choices of $B_x$ and $B_y$, we have
	$d_H(x_i)=2=d_H(y_j)$ and
	$x_iy_j\not\in E(G)$. Then \[d_G(x_i,y_j) =d_H(x_i,y_j) + \sum_{C\in \mathcal{C}}d_C(x_i,y_j) ,\mbox{ where }  i,j\in \{1,2\}.\]
	Since $\delta_2(G)\geqslant 4s$,  we have 
	\[\sum_{C\in \mathcal{C}}d_C(x_i,y_j) =d_G(x_i,y_j) -d_H(x_i,y_j) \geqslant 4(s-1), \mbox{ where }  i,j\in \{1,2\}.\]
	Thus 
	\[\sum_{C\in \mathcal{C}}\sum_{i,j\in \{1,2\}}d_C(x_i,y_j)  \geqslant 16(s-1),\]
	and so either $\sum_{i,j\in \{1,2\}}d_{C'}(x_i,y_j) \geqslant 17$ for some $C'\in \mathcal{C}$,
	or $\sum_{i,j\in \{1,2\}}d_{C}(x_i,y_j)=16$ for all $C \in \mathcal{C}$.

	{\it Case} 1. 
	Suppose first that there is  $C'\in \mathcal{C}$ such that $\sum_{i,j\in \{1,2\}}d_{C'}(x_i,y_j) \geqslant 17$.
	Then, without loss of generality, we let  $ d_{C'}(x_1,y_1) \geqslant 5$.
	By Lemma  \ref{6-cycle},  $C'$ has length $6$. By  (3) of Lemma \ref{degree-3},    $d_{C'}(x_i)\leqslant 3$ and  $d_{C'}(y_j)\leqslant 3$ for all $i,j\in \{1,2\}$.
	In particular, since  $d_{C'}(x_1,y_1) \geqslant 5$,  
	one of $d_{C'}(x_{1})$ and  $d_{C'}(y_{1})$ is $3$, and  the other one is either $2$ or $3$.
	
	Suppose that  $d_{C'}(x_2)=3$.  Write 
	${C'}=u_{1}u_{2}\cdots u_{6}u_1$. 
	By (3) of Lemma \ref{degree-3}, without loss of generality, we may let
	$N_{C'}(x_2)=\{u_1,u_3,u_5\}$.
	Since  $N_{C'}(x_1)\ne\emptyset$,  there is $i\in \{1,3,5\}$ such that $G$ contains a chorded cycle with a chord $x_1x_2$ and vertex set 
	$V(B_x)\cup \{u_i\}$ or   $V(B_x) \cup \{u_i,u_{i+1}\}$, which contradicts Lemma \ref{C-mini}.
	Thus, $d_{C'}(x_2)\leqslant 2$.
	On the other hand, if $d_{C'}(x_1)=3$ and $N_{C'}(x_2)\ne\emptyset$ then we have a similar contradiction.
	
	The argument above says that either $d_{C'}(x_1)=3$ and $d_{C'}(x_2)=0$, or $d_{C'}(x_1)=2$ and $d_{C'}(x_2)\leqslant 2$. Similarly, either $d_{C'}(y_1)=3$ and $d_{C'}(y_2)=0$, or $d_{C'}(y_1)=2$ and $d_{C'}(y_2)\leqslant 2$. If $d_{C'}(x_{1})=d_{C'}(y_{1})=3$ then
	$17\leqslant \sum_{i,j\in \{1,2\}} d_{C'}(x_i,y_j) =d_{C'}(x_1,y_1) + d_{C'}(x_1)+d_{C'}(y_1)\leqslant 6+3+3=12$, a contradiction. If   $d_{C'}(x_{1})=3$ and  $d_{C'}(y_{1})=2$ then
	$17\leqslant \sum_{i,j\in \{1,2\}}d_{C'}(x_i,y_j) =d_{C'}(x_1,y_1) + d_{C'}(x_1,y_2) +d_{C'}(y_1)+d_{C'}(y_2)\leqslant 5+5+2+2=14$, a contradiction. 
	If   $d_{C'}(x_{1})=2$ and  $d_{C'}(y_{1})=3$ then
	$17\leqslant \sum_{i,j\in \{1,2\}}d_{C'}(x_i,y_j) = d_{C'}(x_1,y_1) + d_{C'}(x_2,y_1) +d_{C'}(x_1)+d_{C'}(x_2)\leqslant 5+5+2+2=14$, again a contradiction.

	{\it Case} 2.  Now suppose that $\sum_{i,j\in \{1,2\}}d_{C}(x_i,y_j) = 16$ for all $C\in \mathcal{C}$. 
	If $d_{C'}(x_{i'},y_{j'})  \geqslant 5$ for some $C'\in \mathcal{C}$ and $i',j'\in \{1,2\}$, then an argument similar to that in Case 1 implies that $\sum_{i,j\in \{1,2\}}d_{C'}(x_i,y_j) \leqslant 14$, a contradiction. Thus $d_{C}(x_i,y_j) = 4$ for all $C\in \mathcal{C}$ and $ i,j\in \{1,2\}$.
	We next discuss two cases.

	{\it Subcase} 2.1.
	Suppose that there is $C'\in \mathcal{C}$ such that
	$d_{C'}(x_i) \geqslant 1$ and $d_{C'}(y_j)\geqslant 1$ for some $i,j\in \{1,2\}$.
	Without loss of generality, we  
	choose $C'\in  \mathcal{C}$ with $d_{C'}(x_1)\geqslant 1$, $d_{C'}(y_1)\geqslant 1$ and
	$d_{C'}(x_1)\geqslant d_{C'}(x_2)$.
	Write $C'=u_1u_2\cdots u_tu_1$, and let $u_1x_1\in E(G)$.

	(2.1.1).  Suppose that   $d_{C'}(x_1)=1$. Then $N_{C'}(x_1)=\{u_1\}$, and $d_{C'}(y_j)\geqslant 3$ for all $j\in \{1,2\}$. By Lemma \ref{degree-3}, $t\in \{4,5,6\}$.
	Then either $y_1$ and $y_2$ have a common neighbor on $C'$, or $t=6$ and $y_1$, $y_2$ and the ends of an edge on $C'$ are connected by a $5$-cycle. The former case produces a $4$-cycle with a chord $y_1y_2$, and the latter gives rise to a $5$-cycle with a chord $y_1y_2$.  
	Since $\mathcal{C}$ is minimal, it follows from Lemma \ref{C-mini} that $t=4$.
	
	Recalling that $d_{C'}(x_1,y_j) = 4$, we have $\{u_2,u_3,u_4\}\subseteq N_{C'}(y_j)$, where
	$j\in \{1,2\}$. Then $C_j=u_2y_ju_4u_3u_2$ is a $4$-cycle with a chord $y_ju_3$, where
	$j\in \{1,2\}$.

	Suppose  first that $d_{C'}(y_j)=4$ for some $j\in \{1,2\}$.
	Then
	$N_{C'}(y_j)=\{u_1,u_2,u_3,u_4\}$.
	Pick a shortest path $P[x,y]$ that connects $x$ and $y$ in $H$.
	We have a cycle $D_j=u_1x_1x_2x\cup P[x,y]\cup yy_ju_1$, which has a chord $x_1x$.
	Then $\mathcal{C}\cup\{C_i,D_j\}\setminus\{C'\}$ is a collection of $s$ vertex-disjoint chorded cycles, where $\{i,j\}=\{1,2\}$.
	
	Now let $d_{C'}(y_1)=d_{C'}(y_2)=3$, i.e., $N_{C'}(y_1)=N_{C'}(y_2)=\{u_2,u_3,u_4\}$.
	Since  $d_{C'}(x_2,y_j) = 4$ for $j\in \{1,2\}$, we have $N_{C'}(x_1)=\{u_1\}\subseteq N_{C'}(x_2)$.  Then the $4$-cycle  $C_0=u_1x_1xx_2u_1$ has a  chord  $x_1x_2$, and $\mathcal{C}\cup\{C_0,C_1\}\setminus\{C'\}$ consists of $s$ vertex-disjoint chorded cycles. 
	
	(2.1.2). Suppose that   $d_{C'}(x_1)=2$, and let $N_{C'}(x_1)=\{u_1,u_k\}$.
	Then  $d_{C'}(y_1)\geqslant 2$ and $d_{C'}(y_2)\geqslant 2$.
	In addition, since  $d_{C'}(x_1)\geqslant d_{C'}(x_2)$, we have  $d_{C'}(x_2)\in \{0,1,2\}$.

	Suppose that   $d_{C'}(x_2)=0$. Then $d_{C'}(y_1)=d_{C'}(y_2)=4$, 
	and so $t=4$ by Lemma \ref{degree-3}. We have $N_{C'}(y_1)=N_{C'}(y_2)=\{u_1,u_2,u_3,u_4\}$.
	Then 
	$\mathcal{C}\cup\{C_2,D_1\}\setminus\{C'\}$ is a collection of $s$ vertex-disjoint chorded cycles, where $C_2$ and $D_1$ are constructed as in   (2.1.1).
	
	Suppose first that $d_{C'}(x_2)=1$. 
	Since $d_{C'}(x_2,y_j) = 4$ for all   $ i,j\in \{1,2\}$, we have 
	$d_{C'}(y_1)\geqslant 3$ and $d_{C'}(y_2)\geqslant 3$. 
	By a similar argument as in the first paragraph of   (2.1.1),
	we deduce $t=4$, i.e., $C'$ has length $4$. 
	Since $d_{C'}(x_1,y_1) =d_{C'}(x_1,y_2) = 4$,
	we  have $\{ u_{k_1},u_{k_2}\}\subseteq N_{C'}(y_1)\cap N_{C'}(y_2)$ with
	$\{1,k,k_1,k_2\}=\{1,2,3,4\}$.
	If $N_{C'}(x_2)\subseteq \{u_{k_1},u_{k_2}\}$ then $N_{C'}(y_1)=N_{C'}(y_2)=\{u_1,u_2,u_3,u_4\}$, and
	so  $\mathcal{C}\cup\{C_1,D_2\}\setminus\{C'\}$ is a collection of $s$ vertex-disjoint chorded cycles,
	where $C_1$ and $D_2$ are constructed as  in (2.1.1).
	If $N_{C'}(x_2)\cap \{u_{k_1},u_{k_2}\}=\emptyset$ then, letting  $N_{C'}(x_2)=\{u_1\}$ without loss of generality, 
	$\mathcal{C}\cup\{C_0,C_1\}\setminus\{C'\}$ is a collection of $s$ vertex-disjoint chorded cycles, where $C_0$ and $C_1$ are constructed as  in (2.1.1).

	Next let  $d_{C'}(x_2)=2$. Suppose that $N_{C'}(x_2)\cap N_{C'}(x_1)\ne \emptyset$.
	Without loss of generality, let  $u_1\in N_{C'}(x_2)\cap N_{C'}(x_1)$.
	Then $C_0=u_1x_1xx_2u_1$ is a $4$-cycle with a chord $x_1x_2$. By the choice of $\mathcal{C}$ and Lemma
	\ref{C-mini}, we conclude that $t=4$. Since $d_{C'}(x_1,y_1) =d_{C'}(x_1,y_2) = 4$,
	we  have $\{ u_{k_1},u_{k_2}\}\subseteq N_{C'}(y_1)\cap N_{C'}(y_2)$ with
	$\{1,k,k_1,k_2\}=\{1,2,3,4\}$. Then $u_{k_1}y_1u_{k_2}y_2u_{k_1}$ is a $4$-cycle with a chord $y_1y_2$.
	Thus $C_0$ and  $u_{k_1}y_1u_{k_2}y_2u_{k_1}$ together with $\mathcal{C} \setminus\{C'\}$ form a collection of $s$ vertex-disjoint chorded cycles. Similarly, if  $N_{C'}(y_2)\cap N_{C'}(y_1)\ne \emptyset$ then $G$ contains a collection of $s$ vertex-disjoint chorded cycles.

	Now suppose that $N_{C'}(x_2)\cap N_{C'}(x_1)= \emptyset=N_{C'}(y_2)\cap N_{C'}(y_1)$.
	Write $N_{C'}(x_1,x_2)=\{u_1,u_k,u_{k_1},u_{k_2}\}$ and $N_{C'}(x_1,y_j)=\{u_1,u_k,v_j,w_j\}$, where $j\in \{1,2\}$. Then $u_1, u_k, u_{k_1}, u_{k_2},\\ v_1, v_2, w_1$ and $w_2$ are distinct. In particular, $C'$ has length $t\geqslant 8$. On the other hand, noting that $d_{C'}(x_1,x_2) =4$,
	it is easily checked that there are $z_1\in N_{C'}(x_1)$ and
	$z_2\in N_{C'}(x_2)$ such that $z_1$ and $z_2$ are at distance on $C'$ less than $\frac{t}{2}$.
	Pick the shortest path $P$ on $C'$ that connects $z_1$ and $z_2$.
	Then we get a cycle with a chord $x_1x_2$ and vertex set $V(P)\cup V(B_x)$. By the choice of
	$\mathcal{C}$ and Lemma
	\ref{C-mini}, we have $3+|V(P)|\geqslant t$. Then $3+\frac{t}{2}\geqslant t$, yielding $t\leqslant 6$,
	a contradiction.
	
	(2.1.3). Suppose that   $d_{C'}(x_1)=3$.
	Then $t=|V(C')|  \leqslant 6$ by Lemma \ref{degree-3}. 
	Suppose that $t=5$ or $6$. By   Lemma \ref{degree-3}, $d_{C'}(y_1)\leqslant 3$ and $d_{C'}(y_2)\leqslant 3$, and so $d_{C'}(x_2)\geqslant 1$.
	If $N_{C'}(x_1)\cap N_{C'}(x_2)\ne \emptyset$ then
	there exists a chorded $4$-cycle contained in  $B_x\sqcup u$ with $u\in N_{C'}(x_1)\cap N_{C'}(x_2)$,
	contrary to Lemma \ref{C-mini}. 
	Thus $N_{C'}(x_1)\cap N_{C'}(x_2)= \emptyset$. 
	It follows from (2) and (3) of Lemma \ref{degree-3} that
	there exists an edge on $C'$, say $u_1u_2$ without loss of generality,
	whose ends are adjacent to $x_1$ and $x_2$, respectively.
	Thus we have  a $5$-cycle with a chord $x_1x_2$ and vertex set $\{u_1,u_2,x,x_1,x_2\}$.
	By the choice of $\mathcal{C}$ and Lemma \ref{C-mini}, we have $t=5$.
	In view of (2) of Lemma \ref{degree-3},  we may let $N_{C'}(x_1)=\{u_1,u_3,u_4\}$ and $\{u_2,u_5\} \subseteq N_{C'}(x_2)$. Then $d_{C'}(y_j)\geqslant 2$, and $N_{C'}(y_j) \subseteq \{u_1,u_3,u_4\}$ for $j \in \{1,2\}$.
	Thus $y_1$ and $y_2$ have a common neighbor on $C'$, which yields a $4$-cycle with a chord $y_1y_2$,
	contrary to Lemma \ref{C-mini}. 
	
	The argument above says that $t=4$. Without loss of generality, let $N_{C'}(x_1)=\{u_1,u_2,u_3\}$.
	Then $u_4\in N_{C'}(y_1)\cap N_{C'}(y_2)$, and $u_4y_1yy_2u_4$ is a $4$-cycle with a
	chord $y_1y_2$. In addition, $u_1u_2u_3x_1u_1$ is a $4$-cycle with a
	chord $x_1u_2$. Then these two cycles 
	together with $\mathcal{C} \setminus\{C'\}$ form a collection of $s$ vertex-disjoint chorded cycles.

	(2.1.4). Suppose that   $d_{C'}(x_1)=4$.
	Then $t =|V(C')|=4$ by Lemma \ref{degree-3}. 
	Since $d_{C'}(y_1)\geqslant 1$, without loss of generality, let $u_4\in N_{C'}(y_1)$.
	If $u_4 \in N_{C'}(y_2)$ then two $4$-cycles $yy_1u_4y_2y$ with a chord $y_1y_2$ and $u_1u_2u_3x_1u_1$ with a chord $x_1u_2$ guarantee a collection of $s$ vertex-disjoint chorded cycles. 
	Thus $u_4 \notin N_{C'}(y_2)$ and then $u_4 \in N_{C'}(x_2)$.
	Pick a shortest path $P[x,y]$ that connects $x$ and $y$ in $H$.
	Then we have a cycle $u_4x_2x\cup P[x,y]\cup yy_2y_1u_4$ with a chord $y_1y$,
	and a  $4$-cycle  $u_1u_2u_3x_1u_1$  with a chord $x_1u_2$.
	These two chorded cycles again guarantee a collection of $s$ vertex-disjoint chorded cycles. 
	This completes the proof.

	{\it Subcase} 2.2. 
	Suppose now that for every $C\in \mathcal{C}$, either $d_{C}(x_1)=d_{C}(x_2)=0$ or $d_{C}(y_1)=d_{C}(y_2)=0$.
	In particular, either $d_{C}(y_1)=d_{C}(y_2)=4$ or $d_{C}(x_1)=d_{C}(x_2)=4$, respectively. 
	By Lemma \ref{degree-3}, every $C\in \mathcal{C}$ has length $4$.
	Let $\mathcal{C}_x=\{C\in \mathcal{C}\,:\, d_{C}(y_1)=d_{C}(y_2)=0\}$ and 
	$\mathcal{C}_y=\{C\in \mathcal{C}\,:\, d_{C}(x_1)=d_{C}(x_2)=0\}$.
	Clearly, $\mathcal{C}_x \cap \mathcal{C}_y=\emptyset$ and $\mathcal{C}=\mathcal{C}_x \cup \mathcal{C}_y$.
	If  $\mathcal{C}_x=\emptyset$ or $\mathcal{C}_y=\emptyset$ or $E(V(\mathcal{C}_x),V(\mathcal{C}_y))=\emptyset$, then one of $y$ and $x$ is a cut-vertex of $G$, contrary to the $2$-connectivity of $G$. Thus neither $\mathcal{C}_x$ nor $\mathcal{C}_y$ is empty, and
	we may choose $C_x\in \mathcal{C}_x$ and $C_y\in \mathcal{C}_y$ such that $E(V(C_x),V(C_y))\ne \emptyset$.
	Write $C_x=u_1u_2u_3u_4u_1$, $C_y=v_1v_2v_3v_4v_1$, and let $u_4v_4\in E(G)$.
	We have two $4$-cycles $C_x'=u_1u_2u_3x_1u_1$ and  $C_y'=v_1v_2v_3y_1v_1$, which have chords
	$x_1u_2$ and $y_1v_2$ respectively. Let $\mathcal{C}'=\mathcal{C}\cup\{C_x',C_y'\}\setminus\{C_x,C_y\}$, and $H'=\langle V(H)\cup\{u_4,v_4\}\setminus\{x_1,y_1\}\rangle$.
	Then $|\mathcal{C}'|=|\mathcal{C}|=s-1$, $|V(\mathcal{C}')|=|V(\mathcal{C})|$, $|V(H)|=|V(H')|$ and $H'$ is a 
	component of maximal order in $G-\mathcal{C}'$. 
	Then $r(G,\mathcal{C}')=r(G,\mathcal{C})$, and so $\mathcal{C}'$ is an optimal $(s-1)$-system  of chorded cycles. In addition, it is easy to see that  $H$ has at least one more block than $H'$.

	Suppose that  $H'$  contains two triangle leaf blocks, say $B_{x'}$ and $B_{y'}$, where
	$x'$ and $y'$ are cut vertices of $H'$.   Write $V(B_{x'})=\{x',x'_1,x'_2\}$ and $V(B_{y'})=\{y',y'_1,y'_2\}$. Noting that $\mathcal{C}'$ contains only $4$-cycles,
	we have $\sum_{i,j\in \{1,2\}}d_{C}(x_i',y_j')\leqslant  16$ for all $C\in \mathcal{C}'$. 
	On the other hand, $\sum_{C\in \mathcal{C}'}\sum_{i,j\in \{1,2\}}d_C(x_i',y_j')  \geqslant 16(s-1)$.
	Then $\sum_{i,j\in \{1,2\}}d_{C}(x_i',y_j')=16$ for all $C\in \mathcal{C}'$.  
	If  $d_{C'}(x_i')\geqslant 1$ and $d_{C'}(y_j')\geqslant 1$ for some $i,j\in \{1,2\}$ and some $C'\in \mathcal{C}'$,
	then a similar argument as in Subcase 2.1 implies that $G$ contains a collection of $s$ vertex-disjoint chorded cycles. Thus we may suppose that for every $C\in \mathcal{C}'$, either $d_{C}(x_1')=d_{C}(x_2')=0$ or $d_{C}(y_1')=d_{C}(y_2')=0$. Then, by a similar argument as in the above paragraph, there is an optimal $(s-1)$-system $\mathcal{C}''$ of chorded cycles and a component $H''$ of maximal order in $G-\mathcal{C}''$ such that    $H'$ has at least one more block than $H''$.
	Of course, $|V(H)|=r(G,\mathcal{C})=r(G,\mathcal{C}'')=|V(H'')|$, and    $\mathcal{C}''$ contains only $4$-cycles.

	An inductive repetition of the argument above yields an optimal $(s-1)$-system $\mathcal{C}^*$ and a component $H^*$ of maximal order in $G-\mathcal{C}^*$ such that  $H^*$ has at most one triangle leaf block. Of course, $|V(H)|=r(G,\mathcal{C})=r(G,\mathcal{C}^*)=|V(H^*)|$, and $\mathcal{C}^*$ contains only $4$-cycles.
	For distinct $u,v\in V(H^*)$ with $uv\not\in  E(H^*)$, we have
	$d_{H^*}(u,v)=d_G(u,v)-\sum_{C\in \mathcal{C}^*}d_{C}(u,v)\geqslant 4s-4(s-1)=4$.
	This says that  $\delta_2(H^*)\geqslant 4$. By Corollary \ref{Basic-cor} and Lemma \ref{3-neighbor}, either $H^*$
	contains a chorded cycle, or $H^*$ contains a leaf block that has at least four vertices.
	The former says that  $G$ contains a collection of $s$ vertex-disjoint chorded cycles.
	Suppose that $H^*$ contains a leaf block $B_{x^*}$ with $|V(B_{x^*})|\geqslant 4$, where $x^*$
	is the unique cut-vertex of $H^*$ contained in $B_{x^*}$.
	Then the pair $(B_{x^*},x^*)$ satisfies the hypothesis in Corollary 
	\ref{3-neighbor-cor}, and so  $B_{x^*}$ contains a chorded cycle. 
	Thus  $G$ contains a collection of $s$ vertex-disjoint chorded cycles.
	This completes the proof.
	\qed

	\begin{lemma}\label{Hamiltonian-path}
		Let $G$ be a $2$-connected graph of order at least $4s$ and $\delta_{2}(G) \geqslant 4s$, where $s \geqslant 2$. Let $\mathcal{C}$ be an optimal $(s-1)$-system of chorded cycles in $G$, and let $H$ be a component of maximal order in $G-\mathcal{C}$. 
		Suppose that $H$ contains a Hamiltonian path $P=x_{1}x_{2} x_{3}x_{4}\cdots x_{p}$, where $p\geqslant 5$.  Then $G$ contains a collection of $s$ vertex-disjoint chorded cycles.
	\end{lemma}
	\proof Suppose that the lemma is false. Then  $H$ contains no chorded cycles. This leads to the following observations.
	\begin{claim}\label{claim-HP-1}
		$d_H(x_i)\leqslant 2$ for $i\in \{1,p\}$, $d_H(x_i)\leqslant 3$ for $i\in \{2,p-1\}$, and $d_H(x_i)\leqslant 4$ for $3 \leqslant i \leqslant p-2$.
		For distinct edges $x_ix_j,\, x_{i'}x_{j'}\in E(H)$, if $i<i'<j<j'$ then   $j-i' \geqslant 2$.
	\end{claim}
	In view of Claim \ref{claim-HP-1},   if  $x_1x_3$ or $x_{p-2}x_p\in H$ then 
	it is easy to see that $x_1x_2x_3x_1$ or $x_{p-2}x_{p-1}x_px_{p-2}$ is a leaf block of $H$.
	Then  Lemma \ref{triangle-leaf-block} implies the following assertion.
	\begin{claim}\label{claim-HP-2}
		One of $x_1x_3$ and $x_{p-2}x_p$,  say $x_1x_3$ without loss of generality,  is not an edge of $H$.
	\end{claim}

	\begin{claim}\label{claim-HP-3} There exist no consecutive vertices $x_i$'s 
		such that one of the following holds:
		\begin{itemize}
			\item[(1)] $x_1x_3,\,x_1x_4,\,x_2x_4\not\in E(H)$, and $d_H(x_1,x_3)+d_H(x_1,x_4)+d_H(x_2,x_4)\leqslant 11$;
			\item[(2)] $x_1x_3,\,x_2x_4,\,x_3x_5\not\in E(H)$, and $d_H(x_1,x_3)+d_H(x_2,x_4)+d_H(x_3,x_5)\leqslant 11$;
			\item[(3)] $x_1x_3,\,x_1x_4,\,x_2x_5\not\in E(H)$ and $x_2x_4\in E(H)$, and $d_H(x_1,x_3)+d_H(x_1,x_4)+d_H(x_2,x_5)\leqslant 11$.
		\end{itemize}
	\end{claim}
	\noindent{\it Proof of Claim} \ref{claim-HP-3}.
	Suppose, to the contrary, that there exist consecutive vertices $x_i$'s on $P$ that
	satisfy one of (1)-(3) above. For convenience, denote by $(x_{i_1},x_{j_1})$, 
	$(x_{i_2},x_{j_2})$  and  $(x_{i_3},x_{j_3})$  the three pairs of nonadjacent vertices in  (1), (2) or (3).  Then, since $\delta_2(G)\geqslant 4s$, we have
	\[\sum_{C\in \mathcal{C}}\sum_{a=1}^3d_C(x_{i_a},x_{j_a})=
	\sum_{a=1}^3d_G(x_{i_a},x_{j_a})-\sum_{a=1}^3d_H(x_{i_a},x_{j_a})
	\geqslant 12s-11=12(s-1)+1.\]
	We get $\sum_{a=1}^3d_C(x_{i_a},x_{j_a})\geqslant 13$ for at least one chorded cycle $C\in \mathcal{C}$, which contradicts Lemma \ref{5-path}. Thus Claim \ref{claim-HP-3} follows.
	\qed
	
	Based on the claims above, we next deduce a contradiction.
	By Claims \ref{claim-HP-1} and \ref{claim-HP-2},  $d_H(x_1)\leqslant 2$,  $d_H(x_2)\leqslant 3$ 
	and $d_H(x_3)\leqslant 3$. We shall get the  contradiction  in  two cases, say  $d_H(x_1)=1$, and
	$d_H(x_1)=2$.

	{\it Case }1. Suppose that $d_H(x_1)=1$. 
	If $d_H(x_2)=2$ then  $x_1x_3,\,x_1x_4,\, x_2x_4\not\in E(H)$, 
	$d_H(x_1,x_3) \leqslant 3$, $d_H(x_1,x_4)\leqslant 4$ and $d_H(x_2,x_4)\leqslant 4$, 
	and so $d_H(x_1,x_3)+d_H(x_1,x_4)+d_H(x_2,x_4)\leqslant 11$, contrary to (1) of Claim \ref{claim-HP-3}.
	If $x_2x_4\in E(H)$ then $x_1x_3,\,x_1x_4,\, x_2x_5\not\in E(H)$,  
	$d_H(x_1,x_3)=2$, $d_H(x_1,x_4)\leqslant 4$ and $d_H(x_2,x_5)\leqslant 5$, 
	and so $d_H(x_1,x_3)+d_H(x_1,x_4)+d_H(x_2,x_5)\leqslant 11$, contrary to (3) of Claim \ref{claim-HP-3}.
	Thus we suppose further that $d_H(x_2)=3$ and $x_2x_4\not\in E(H)$.
	Write $N_H(x_2)=\{x_1,x_3,x_k\}$ for some $k\geqslant 5$. 
	
	Using Claim \ref{claim-HP-1}, it is easily observed that
	$d_H(x_i)\leqslant 3$ for all $3\leqslant i \leqslant k-1$, and the equality holds for at most one
	$i$.   Also, if $d_H(x_i)=3$ for   $3\leqslant i \leqslant k-1$, then $3\leqslant i \leqslant k-2$ and
	$x_i$ has a neighbor lying  on the path $P[x_{k+1},x_p]$. 
	If  such an $i$ exists then denote it by $i_0$, and put $i_0=2$ otherwise. 
	We have $x_1x_3,\,x_1x_4,\, x_2x_4\not\in E(H)$, and 
	\[ \left\{
	\begin{array}{cl}
		d_H(x_1,x_3)+d_H(x_1,x_4)+d_H(x_2,x_4)\leqslant  2+3+4=9  & \mbox{ if   }i_0=2   \mbox{ or } i_0 \geqslant 5,   \\
		d_H(x_1,x_3)+d_H(x_1,x_4)+d_H(x_2,x_4)\leqslant   3+3+4=10  &  \mbox{ if  }i_0=3,  \\
		d_H(x_1,x_3)+d_H(x_1,x_4)+d_H(x_2,x_4)\leqslant   2+4+5=11 &    \mbox{ if  }i_0=4. 
	\end{array}
	\right.\]
	Clearly, each case leads to a contradiction.

	{\it Case }2. Suppose that $d_H(x_1)=2$.  Write $N_H(x_1)=\{x_2,x_k\}$ for some $k\geqslant 4$.

	Suppose that $k=4$. By Lemma \ref{v2}, $d_{H}(x_{3})=2$, and $x_{1}x_{3},\,x_{2}x_{4},\,x_{3}x_{5} \notin E(H)$. It follows that $d_{H}(x_{1},x_{3})+ d_{H}(x_{2},x_{4}) + d_{H}(x_{3},x_{5}) \leqslant 2+5+4=11$, contrary to (2) of Claim \ref{claim-HP-3}.
	Now let $k\geqslant 5$.
	Claim \ref{claim-HP-1} leads to an observation similar to that made as  Case 1,
	that is, $d_H(x_i)=3$ holds for at most one $i$ from $2$ to $k-2$, and $d_H(x_i)=2$ for any other
	$i$ from $1$ to $k-1$, which in turn implies that $x_1x_3,\,x_{1}x_{4},\,x_2x_4\not\in E(H)$.
	If $d_{H}(x_{4})=2$ then it is easily checked that $d_{H}(x_{1},x_{3})+d_{H}(x_{1},x_{4})+d_{H}(x_{2},x_{4}) \leqslant 11$, contrary to (1) of Claim \ref{claim-HP-3}. 
	Thus $d_{H}(x_{4})=3$. We have $k \geqslant 6$, and $d_{H}(x_{i})=2$ for $i\in \{2,3,5\}$. It follows that $d_{H}(x_{1},x_{3})+ d_{H}(x_{2},x_{4}) + d_{H}(x_{3},x_{5}) \leqslant 3+4+3=10$, contrary to (2) of Claim \ref{claim-HP-3}.
	This completes the proof.
	\qed


	\vskip 20pt
	
	\section{The proof of Theorem \ref{T4}}

	Suppose that  Theorem \ref{T4} is false,  and let $G$ be a counterexample of minimal order.
	Clearly, $G$ is not a complete graph and, in view of Corollary \ref{Basic-cor},
	$s\geqslant 2$. Further, we choose $G$ with $|E(G)|$ as large as possible.
	Pick two nonadjacent vertices $x$ and $y$ in $G$, and
	let $G_{xy}$ be the graph obtained from $G$ by adding an edge that joins the chosen vertices $x$ and $y$.
	Then, by the choice of $G$, there exists an optimal $s$-system $\mathcal{C}_{xy}$ of 
	chorded cycles in $G_{xy}$. Since $G$  is a counterexample, 
	$x$ and $y$ lie on the same cycle $C_{xy}$ of these $s$-cycles.
	Thus we  have  a collection $\mathcal{D}=\mathcal{C}_{xy}\setminus \{C_{xy}\}$ of $s-1$
	vertex-disjoint chorded cycles in $G$, and $|V(G-\mathcal{D})|\geqslant 4$.
	This implies that $G$ contains an  optimal $(s-1)$-system $\mathcal{C}$ of 
	chorded cycles such that $|V(G-\mathcal{C})|\geqslant 4$.
	
	\vskip 5pt
	
	In the following, 
	we let  $\mathcal{C}$ be an  optimal $(s-1)$-system  of 
	chorded cycles in $G$ with $|V(G-\mathcal{C})|\geqslant 4$, and let
	$H$ be a component of maximal order in $G-\mathcal{C}$. In particular,  $r(G,\mathcal{C})=|V(H)|$. Clearly, $H$ does not contain chorded cycles, and two vertices of $H$ are adjacent in $H$ if and only if 
	they are adjacent in $G$. 
	
	\begin{claim}\label{comp-H-1}
		If $u\in H$, $v\in V(G-\mathcal{C})\setminus V(H)$ then $d_C(u,v)\leqslant 4$ for all $ C\in \mathcal{C}$. 
	\end{claim}
	\noindent{\it Proof of Claim} \ref{comp-H-1}. Suppose not, and let $ C\in \mathcal{C}$ with $d_C(u,v)\geqslant 5$.
	By Lemma \ref{6-cycle}, $|V(C)|=6$, and $C_v=(C-u')\sqcup v$ is a chorded $6$-cycle, where $u' \in N_{C}(u)$.
	Replacing $C$ with $C_v$, we have a minimal $(s-1)$-system $\mathcal{C}_v$ in $G$. However, $G-\mathcal{C}_v$ has a component with vertex set $V(H)\cup\{u'\}$, 
	and so $r(G,\mathcal{C})< r(G,\mathcal{C}_v)$,
	contrary to the optimality of $\mathcal{C}$.  Hence the claim is proven.
	\qed

	\begin{claim}\label{comp-H-2}
		$|V(H)|\geqslant 3$.
	\end{claim}
	\noindent{\it Proof of Claim} \ref{comp-H-2}.
	Suppose  that  $|V(H)|=1$. Picking distinct $u,v\in  V(G-\mathcal{C})$, we have $uv\not\in E(G)$, and 
	$\sum_{C\in \mathcal{C}} d_C(u,v)=d_G(u,v)\geqslant 4s=4(s-1)+4$. Then there is $C\in \mathcal{C}$ such that $d_C(u,v)\geqslant 5$; however, $d_C(u,v)\leqslant 4$ by Claim \ref{comp-H-1}, a contradiction.

	Suppose  that  $|V(H)|=2$. Since $|V(G-\mathcal{C})|\geqslant 4$, pick 
	a   component $K$ in $G-\mathcal{C}$ other than $H$. Then $|V(K)|\leqslant 2$,  
	$uv\notin E(G)$ for all $u\in V(H)$ and $v\in V(K)$, and $\sum_{u\in V(H),v\in V(K)}d_H(u,v)=
	2|V(K)|^2$. Write $V(H)=\{u,u'\}$. We have 
	\[\sum_{C\in \mathcal{C}} \sum_{v\in V(K)}(d_C(u,v)+d_C(u',v))\geqslant 8|V(K)|s-2|V(K)|^2>8|V(K)|(s-1).\]
	Then there is $C\in \mathcal{C}$ such that $\sum_{v\in V(K)}(d_C(u,v)+d_C(u',v))\geqslant  8|V(K)|+1$.
	This in turn implies that there is  $v\in V(K)$ such that $d_C(u,v)+d_C(u',v)\geqslant 9$.
	We have either $d_C(u,v) \geqslant 5$ or $d_C(u',v) \geqslant 5$, which contradicts our Claim \ref{comp-H-1}.
	\qed
	
	\begin{claim}\label{comp-H-3}
		$|V(H)|\geqslant 4$.
	\end{claim}
	\noindent{\it Proof of Claim} \ref{comp-H-3}.
	Suppose  the contrary, then $|V(H)|=3$ by Claim \ref{comp-H-2}. Clearly, $H$ contains a $3$-path, say $u_1u_2u_3$.
	Picking  $v\in  V(G-\mathcal{C})\setminus V(H)$ with degree as small as possible, we have $d_{G-\mathcal{C}}(v)\leqslant 2$, and
	$u_iv\not\in E(G)$ for $i\in \{1,2,3\}$. It is easy to see that
	$d_{G-\mathcal{C}}(u_i,v)\leqslant 4$, and so
	$\sum_{i=1}^3d_{G-\mathcal{C}}(u_i,v)\leqslant 12$.
	If $\sum_{i=1}^3d_{G-\mathcal{C}}(u_i,v)\leqslant 11$ then $\sum_{C\in \mathcal{C}} \sum_{i=1}^3d_C(u_i,v)\geqslant 12s-11=12(s-1)+1$. 
	Thus there is $C'\in \mathcal{C}$ such that $\sum_{i=1}^3d_{C'}(u_i,v)\geqslant 13$, and  one of the three summands is at least $5$. 
	Combining this with Claim \ref{comp-H-1}, we obtain a contradiction.
	This forces that $\sum_{i=1}^3d_{G-\mathcal{C}}(u_i,v) =12$ and $d_{G-\mathcal{C}}(u_i,v) =4$. It follows that $H$ is a triangle and $d_{G-\mathcal{C}}(v)=2$. Let $K$ be the component of $G-\mathcal{C}$ containing $v$. Then $|V(K)|\leqslant |V(H)|=3$. Recalling that $v$ has degree as small as possible in $G-\mathcal{C}$ and $d_{G-\mathcal{C}}(v)=2$,
	we have $|V(K)|=3$, and further $K$ is  a triangle. Let $V(K)=\{v_{1},v_{2},v_{3}\}$.
	
	Fix $j \in \{1,2,3\}$. By Claim \ref{comp-H-1}, we have $\sum_{i=1}^3d_{C}(u_i,v_j) \leqslant  12$. 
	Then $12s \leqslant \sum_{i=1}^3d_{G}(u_i,v_j) =  \sum_{i=1}^3d_{G-\mathcal{C}}(u_i,v_j)+\sum_{C\in \mathcal{C}}\sum_{i=1}^3d_{C}(u_i,v_j)\leqslant12+12(s-1)=12s$. Thus $\sum_{i=1}^3d_{C}(u_i,v_j) =12$ and $d_{C}(u_i,v_j)=4$ for any $C \in \mathcal{C}$.

	{\it Case }1. Suppose first that for every $C\in \mathcal{C}$, either $d_{C}(u_1)=d_{C}(u_2)=d_{C}(u_3)=0$ or $d_{C}(v_1)=d_{C}(v_2)=d_{C}(v_3)=0$.
	In particular, either $d_{C}(v_1)=d_{C}(v_2)=d_{C}(v_3)=4$ or $d_{C}(u_1)=d_{C}(u_2)=d_{C}(u_3)=4$, respectively. 
	By Lemma \ref{degree-3}, every $C\in \mathcal{C}$ has length $4$.
	Let $\mathcal{C}_u=\{C\in \mathcal{C}\,:\, d_{C}(v_1)=d_{C}(v_2)=d_{C}(v_3)=0\}$ and 
	$\mathcal{C}_v=\{C\in \mathcal{C}\,:\, d_{C}(u_1)=d_{C}(u_2)=d_{C}(u_3)=0\}$.
	Clearly, $\mathcal{C}_u \cap \mathcal{C}_v=\emptyset$ and $\mathcal{C}=\mathcal{C}_u \cup \mathcal{C}_v$.
	Let $M(V(\mathcal{C}_u),V(\mathcal{C}_v))$ denote a matching between $\mathcal{C}_u$ and $\mathcal{C}_v$.
	If  $\mathcal{C}_u=\emptyset$ or $\mathcal{C}_v=\emptyset$ or $|M(V(\mathcal{C}_u),V(\mathcal{C}_v))|\leqslant 1$, then $G$ is not $2$-connected, a contradiction. Thus neither $\mathcal{C}_u$ nor $\mathcal{C}_v$ is empty, and $|M(V(\mathcal{C}_u),V(\mathcal{C}_v))|\geqslant 2$.
	Note that $G-\mathcal{C}$ has only two components, otherwise it is readily verified that $G$ contains $s$ vertex-disjoint chorded cycles, a contradiction. Thus $|V(G)|=4s+2$.
	By the neighborhood-union condition $\delta_{2}(G) \geqslant 4s$, we have that 
	\begin{equation*}
		G \cong K_{4\ell+3} \cup K_{4(s-1-\ell)+3} \cup M,~~for~~ \ell \in \{1,\cdots, s-2\}.
	\end{equation*}
	where $M$ is a matching and $|M| \geqslant 2$.
	Thus $G$ has $s$ vertex-disjoint chorded cycles, a contradiction.

	{\it Case }2. Then suppose that there is $C\in \mathcal{C}$ such that $d_{C}(u) \geqslant 1$ and $d_{C}(v) \geqslant 1$ for some $u \in V(H)$ and $v \in V(K)$. We assume that $d_{C}(u_{1}) \geqslant d_{C}(u_{2})  \geqslant d_{C}(u_{3})$. Let $C=w_{1}w_{2} \cdots w_{t}w_1$.
	
	Suppose that there are $u,u' \in V(H)$ such that $d_{C}(u)=d_{C}(u')=3$. Then $|V(C)| \leqslant 6$  by Lemma \ref{degree-3}. If $|V(C)|=5$ then $u$ and $u'$ share a  neighbor on $C$, and so we obtain a chorded $4$-cycle, contrary to Lemma \ref{C-mini}.
	If $|V(C)|=6$ then either $u$ and $u'$ share a  neighbor on $C$ or they are respectively adjacent to the ends of an edge of $C$, each of two cases gives rise to a chorded  cycle
	of length less than $6$, contrary to Lemma \ref{C-mini}.
	Thus $|V(C)| =4$ and $d_{C}(v_{j}) \geqslant 1$ for all $j \in \{1,2,3\}$. Let $N_{C}(u)=\{w_{1},w_{2},w_{3}\}$. Then $w_{4} \in N_{C}(v_{j})$ for all $j \in \{1,2,3\}$. It follows that $C \cup H \cup K$ contains two chorded $4$-cycles $uw_{1}w_{2}w_{3}u$ with chord $uw_{2}$ and $v_{1}v_{3}v_{2}w_{4}v_{1}$ with chord $v_{1}v_{2}$. Now we have $s$ vertex-disjoint chorded cycles, a contradiction.
	Similarly, if there are $v,v' \in V(K)$ such that $d_{C}(v)=d_{C}(v')=3$ then we have a similar contradiction.
	
	Thus there is at most a vertex $u \in V(H)$ and a vertex $v \in V(K)$ such that $d_{C}(u)=d_{C}(v)=3$.
	This implies that $d_{C}(u) \geqslant 2$ for all $u \in V(H)$ and $d_{C}(v) \geqslant 2$ for all $v \in V(K)$.
	
	Suppose that there is $u \in V(H)$ such that $d_{C}(u) =2$. Then  $d_{C}(v_{j}) \geqslant 2$ for all $j \in \{1,2,3\}$.
	If  $C$ has length $5$ then there are distinct $v_i$ and $v_j$
	that have a common neighbor on $C$, and a chorded $4$-cycle arises, contrary to Lemma \ref{C-mini}.
	If $C$ has length $6$ then there are distinct $v_i$ and $v_j$ such that they share a neighbor on $C$, or they are respectively adjacent to the ends of an edge or a chord of $C$, each of the three cases gives rise to a chorded  cycle of length less than $6$,  a contradiction.
	If $C$ has length at least $7$, then there are distinct $v_i$ and $v_j$ such that they 
	are respectively adjacent  with two vertices on $C$ that have distance less than $\frac{|V(C)|}{2}$ on $C$, and thus a chorded cycle of length no more than $|V(C)|-1$ arises, again a contradiction.
	Thus $C$ is a $4$-cycle. Then two of $u_1$, $u_2$ and $u_3$ share a neighbor, say $w$, on $C$, and so $H\sqcup w$ contains a chorded $4$-cycle $C'$. Now we have a minimal $(s-1)$-system 
	$\mathcal{C}'=\mathcal{C}\cup\{C'\}\setminus\{C\}$  of chorded cycles.
	However,  $G-\mathcal{C}'$ has a component with vertex set $V(K)\cup (V(C)\setminus\{w\})$,
	contrary to the optimality of $\mathcal{C}$. 
	Similarly, if there is $v \in V(K)$ such that $d_{C}(v)=2$, then we have a similar contradiction.
	Thus, we may assume that $d_{C}(u) \geqslant 3$ and $d_{C}(v) \geqslant 3$ for all $u \in V(H)$ and $v \in V(K)$.
	
	Recall that there is at most one vertex $u \in V(H)$ with  $d_{C}(u)=3$ and at most one vertex $v \in V(K)$ such that $d_{C}(v)=3$.
	Thus suppose that $d_{C}(u_{2}) =d_{C}(u_{3})=d_{C}(v_{2})=d_{C}(v_{3})=4$. Then, by Lemma \ref{degree-3}, $|V(C)|=4$. It is easy to check that $C \cup H\cup K$ has two vertex-disjoint chorded cycles, a contradiction.
	\qed

	\begin{claim}\label{comp-H-4}
		\begin{itemize}
			\item[(1)] $H$ contains no Hamiltonian paths.
			\item[(2)] There are nonadjacent vertices $u,v\in V(H)$ such that $d_H(u,v)\leqslant 3$, in particular, $\mathcal{C}$ contains a chorded $6$-cycle.
		\end{itemize}
	\end{claim}   
	\noindent{\it Proof of Claim} \ref{comp-H-4}. 
	Since $G$ is a counterexample,  (1) of the claim follows from Lemma \ref{4-path-cycle} and Lemma \ref{Hamiltonian-path}, and the first part of (2)
	follows from Lemmas \ref{3-neighbor} and \ref{triangle-leaf-block}. 
	Pick distinct $u,v\in V(H)$ with $uv\not\in E(G)$ and $d_H(u,v)\leqslant 3$.
	Then $\sum_{C\in \mathcal{C}}d_C(u,v)\geqslant 4s-3$, and so
	$d_C(u,v)\geqslant 5$ for some $C\in \mathcal{C}$. 
	Thus the second part of (2) follows from Lemma \ref{6-cycle}. 
	\qed
	
	\vskip 10pt

	For convenience, for an optimal $(s-1)$-system $\mathcal{C}$ of chorded cycles in $G$, 
	let  $\ell(\mathcal{C})$ be the maximum length of paths contained in the components
	with maximal order in $G-\mathcal{C}$.
	Choose  an optimal $(s-1)$-system $\mathcal{C}$  and a component $H$ with  maximal order in $G-\mathcal{C}$ that satisfy the following assumption:
	\begin{itemize}
		\item[$(\ddag)$]  $\ell(\mathcal{C})$ is as large as possible, and  $H$ contains an $\ell(\mathcal{C})$-path $P=u_1u_2\cdots u_p$, where $p=\ell(\mathcal{C})+1$.
	\end{itemize}
	
	Clearly,  $p\geqslant 3$ as $|V(H)|\geqslant 4$ by Claim \ref{comp-H-3}, and
	$V(P)\ne V(H)$ by (1) of Claim \ref{comp-H-4}. 
	By the choice of $P$, it is easily shown that neither $u_1$ nor $u_p$ is a cut-vertex of $H$. 
	In addition,  $u_1u_p\not\in E(H)$.   Suppose the contrary, then since $H$ is connected, there are $v\in V(H-P)$ and $u_k$ with $2 \leqslant k \leqslant p-1 $ such that $vu_k\in E(H)$. This leads to a $(p+1)$-path 
	$vu_ku_{k-1}\cdots u_1u_pu_{p-1}\cdots u_{k+1}$ in $H$, contrary to the choice of $P$.
	Therefore, $u_1u_p\not\in E(H)$.

	\begin{claim}\label{comp-H-5}
		Suppose that $\mathcal{C}$, $H$ and $P$ satisfy the assumption $(\ddag)$, and there exists a vertex $w \in V(H) \setminus (V(P))$ such that  $d_{H}(w) \leqslant 2$, $u_{1}w, u_{p}w \notin E(H)$ and $H-w$ is connected. Then the following statements hold:
		\begin{itemize}
			\item[$(1)$]   $d_{H}(w) =d_{H}(u_{1})=d_{H}(u_{p})=2$;
			\item[$(2)$]   $d_{C}(u_1,w)=d_{C}(u_p,w)=d_{C}(u_{1},u_{p})=4$  for any $C \in \mathcal{C}$;
			\item[$(3)$]  $d_{C}(w)=2$ for any chorded $6$-cycle $C \in \mathcal{C}$.   
		\end{itemize}
	\end{claim}   
	\noindent{\it Proof of Claim} \ref{comp-H-5}.  
	Suppose that there is a chorded cycle $C' \in \mathcal{C}$ such that $d_{C'}(u_{i},w) \geqslant 5$ for some $i \in \{1,p\}$.  By Lemma  \ref{6-cycle}, $|V(C')|=6$ and there exists $u'\in N_{C'}(u_i)$ such that
	$C_w=(C'-u')\sqcup w$ is a chorded cycle of length $6$. Since $H-w$ is connected, replacing $C'$ by $C_w$, we have an optimal $(s-1)$-system of chorded cycles, say $\mathcal{C}_{w}$. Noting that $u_{i}u' \in E(G)$,  as a component in $G-\mathcal{C}_{w}$, the subgraph $(H-w)\sqcup u'$ contains a path of length $p=\ell(\mathcal{C})+1$, which contradicts $(\ddag)$. 
	Therefore, 
	\begin{equation}\label{eq-1}
		d_{C}(u_{i},w) \leqslant 4,\,\, for \ any \ i \in \{1,p\} \ and \ any \ C \in \mathcal{C}.
	\end{equation}
	
	By (2) of Lemma \ref{degree-2}, $d_H(u_i)=d_P(u_i)\leqslant 2$ for each $i \in \{1,p\}$.
	Since $d_{H}(u_{i},w)= d_{G}(u_{i},w)-\sum_{C\in \mathcal{C}}d_{C}(u_{i},w) \geqslant 4s-4(s-1) \geqslant 4$, we get $d_H(w) \geqslant 2$. By the assumption, $d_{H}(w) =2$, and so $d_{H}(u_{1})=d_{H}(u_{p})=2$. 
	Then (1) of the claim follows.
	
	If $d_{C'}(u_{i},w) \leqslant 3$ for some $i \in \{1,p\}$ and some $C' \in \mathcal{C}$, then $d_{H}(u_{i},w)= d_{G}(u_{i},w)-\sum_{C\in \mathcal{C}\setminus\{C'\}}d_{C}(u_{i},w)-d_{C'}(u_{i},w) \geqslant 4s-4(s-2)-3 \geqslant 5$, which contradicts (1). Thus, by \eqref{eq-1}, $d_{C}(u_{1},w) =d_{C}(u_{p},w) =4$ for any $C\in \mathcal{C}$, desired as in (2) of the claim.

	Since $u_{1}u_{p} \notin E(G)$, we have $\sum_{C\in \mathcal{C}}d_{C}(u_{1},u_{p})\geqslant 4s-d_{H}(u_{1},u_{p})\geqslant 4(s-1)$. Then either $d_{C}(u_{1},u_{p})=4$ for any  $C\in \mathcal{C}$,
	or $d_{C'}(u_{1},u_{p}) \geqslant 5$
	for some $C'\in \mathcal{C}$. 
	Suppose the latter case occurs.  
	Since $u_{1}u_{p} \notin E(G)$,
	by Lemma  \ref{6-cycle}, $|V(C')|=6$, and 
	$d_{C'}(u_{1})=3$ or $d_{C'}(u_{p}) =3$.   Writing  $C'=w_1w_2w_3w_4w_5w_6w_1$,
	without loss of generality, let 
	$d_{C'}(u_{1})= 3$,  $N_{C'}(u_{1})=\{w_{1},w_{3},w_{5}\}$ and $\{w_{2},w_{4}\} \subseteq N_{C'}(u_{p})\subseteq \{w_{2},w_{4},w_6\}$.
	Since $d_{C'}(u_{1},w)=4$, we have $N_{C'}(w)\cap \{w_{2},w_{4},w_6\}\ne\emptyset$.
	Suppose that $w_2\in  N_{C'}( w)$. Then, since   $d_{C'}(u_{p},w)=4$, we have
	$ N_{C'}( w)\supseteq \{w_2, w_j\}$ for some $j\in \{1,3,5\}$.
	Since $C'\sqcup w$ is triangle-free by (3) of
	Lemma \ref{degree-3}, we have $w_j=w_5$ and  $ N_{C'}( w)\supseteq\{w_2, w_5\}$. 
	Again by (3) of Lemma \ref{degree-3}, $w_2w_5\not\in E(G)$  and, since
	$C'$ is a chorded cycle, either $w_1w_4$ or $w_3w_6$ is an edge.
	It follows that $(C'-w_2)\sqcup u_1$ is a chorded $6$-cycle.
	Let $\mathcal{C}_{u_1}=\mathcal{C}\cup\{(C'-w_2)\sqcup u_1\}\setminus\{C'\}$.
	Then $\mathcal{C}_{u_1}$ is a minimal $(s-1)$-system of chorded cycles.
	Recall that $u_1$ is not a cut vertex of $H$.
	It follows that $(H-u_1)\sqcup w_2$ is a component of $G-\mathcal{C}_{u_1}$.
	However, $(H-u_1)\sqcup w_2$ contains a path $u_2\cdots u_pw_2w$ of length $p=\ell(\mathcal{C})+1$,
	which contradicts the assumption $(\ddag)$.
	Similarly, a contradiction arises from either
	$w_4\in  N_{C'}( w)$ or $w_6\in N_{C'}( w)\cap N_{C'}(u_{p})$.
	Thus, $w_6\in N_{C'}( w)\subseteq\{w_1,w_3,w_5,w_6\}$ and $N_{C'}(u_{p})=\{w_{2},w_{4}\}$.
	Since $C'\sqcup w$ is triangle-free by (3) of
	Lemma~\ref{degree-3}, we have $w_1, w_5\not\in N_{C'}( w)$, and so
	$ N_{C'}( w)=\{w_3,w_6\}$. Now,  $w_3w_6\not\in E(G)$ and either $w_1w_4$ or $w_2w_5$ is an edge.
	Then    $(C'-w_3)\sqcup u_p$ is a chorded $6$-cycle, which leads to a similar contradiction. 
	Thus  $d_{C}(u_{1},u_{p})=4$ for any  $C\in \mathcal{C}$, and  (2) of the claim follows.

	Let    $C' \in \mathcal{C}$ be a $6$-cycle,  and write $C'=w_1w_2w_3w_4w_5w_6w_1$.
	By  Lemma \ref{degree-3}, $d_{C'}(u_{1}) \leqslant 3$ and $d_{C'}(w) \leqslant 3$, and so
	$d_{C'}(w) \geqslant 1$ and $d_{C'}(u_{1}) \geqslant 1$ by (2) of the claim.
	If $d_{C'}(w)=1$ then $d_{C'}(u_{1})=d_{C'}(u_{p})=3$,
	which yields that   either $N_{C'}(u_{1})=N_{C'}(u_{p})$ or $N_{C'}(u_{1}) \cup N_{C'}(u_{p})=V(C')$, contrary to $d_{C'}(u_{1},u_{p})=4$.
	Thus $d_{C'}(w)=2$ or $3$.
	Suppose that $d_{C'}(w)=3$. Without loss of generality, let  $N_{C'}(w)=\{w_{1},w_{3},w_{5}\}$ and  $w_{2} \in N_{C'}(u_{p})$. Then  we have a chorded $6$-cycle $ww_1w_6w_5w_4w_3w$, say $C''$, which has a chord $ww_5$. 
	Recalling that $H-w$ is connected, we get an optimal $(s-1)$-system of chorded cycles, say $\mathcal{C}''$, by replacing $C'$ with $C''$. Now  $G-\mathcal{C}''$ has a component $(H-w)\sqcup w_2$ which contains a path of length $p=\ell(\mathcal{C})+1$, contrary to $(\ddag)$. Thus $d_{C'}(w)=2$, and (3) of the claim follows.
	\qed
	
	\vskip 10pt

	From now on, let $\mathcal{C}$, $H$ and $P$ satisfy the assumption $(\ddag)$, and  choose a longest path $Q=v_1v_2\cdots v_q$ in $H-P$ with $d_{P_{1}}(v_1) \leqslant d_{P_{1}}(v_q)$.
	By Lemma \ref{degree-2} and the choices of $P$ and $Q$, for $i \in \{1,p\}$ and $j \in \{1,q\}$, we have $u_iv_j \notin E(H)$, and 
	\begin{equation}\label{eq-2}
		d_{H}(u_i)=d_{P}(u_i) \leqslant 2.
	\end{equation}
	\begin{equation}\label{eq-3}
		d_H(v_j)=d_{P\sqcup Q}(v_j) = d_P(v_j)+d_Q(v_j),\,\, d_{Q}(v_j) \leqslant 2,\,\, d_P(v_1)+d_P(v_q) \leqslant 3.
	\end{equation}
	
	We claim that $d_{H}(v_{1}) \leqslant 2$. Suppose the contrary, then $d_{P}(v_1) \geqslant 1$ by \eqref{eq-3}.
	Recalling  that $d_{P}(v_{1}) \leqslant d_{P}(v_{q})$, we have $d_{P}(v_{q}) \geqslant 1$. 
	If $d_{Q}(v_{1}) \geqslant 2$ or $d_{Q}(v_{q}) \geqslant 2$, then $P \sqcup Q$ has a chorded cycle, a contradiction.
	Thus $d_{Q}(v_{j})=1$, and it follows that $d_{P}(v_{j}) \geqslant 2$ for each $j \in \{1,q\}$.
	Considering two vertices $v_{1}$ and $v_{q}$, each sends at least two edges to another path $P$, and by Lemma \ref{vertex-disjoint-paths}, a chorded cycle exists in $P \sqcup Q$, again a contradiction. 
	Therefore, $d_{H}(v_{1}) \leqslant 2$.  
	
	By \eqref{eq-3}, removing $v_1$ only affects connectivity of $P\sqcup Q $. If $d_{P}(v_1)=0$ then $d_{H}(v_1)=d_{Q}(v_1)$ and $H-v_1$ is connected. For $d_{P}(v_1) \geqslant 1$, recalling that  $d_{P}(v_1) \leqslant d_{P}(v_q)$, we get $d_{P_{1}}(v_q) \geqslant 1$, and so $H-v_1$  is also connected. 
	Thus,   $v_{1}$ is not a cut-vertex of $H$. 
	Now $v_1$ satisfies the hypotheses in Claim \ref{comp-H-5}. Then
	\begin{itemize}
		\item[(1')]   $d_{H}(v_1) =d_{H}(u_{1})=d_{H}(u_{p})=2$; and 
		\item[(2')]   $d_{C}(u_1,v_1)=d_{C}(u_p,v_1)=d_{C}(u_{1},u_{p})=4$  for any $C \in \mathcal{C}$; and 
		\item[(3')]  $d_{C}(v_1)=2$ for any chorded $6$-cycle $C \in \mathcal{C}$.   
	\end{itemize}
	
	\vskip 5pt
	
	Now we are ready to finish  the proof of Theorem \ref{T4} by deriving a final contradiction. 
	By Claim \ref{comp-H-4}, $\mathcal{C}$ contains a chorded $6$-cycle, say $C^{*}= w_1w_2w_3w_4w_5w_6w_1$.
	It follows from   (2') and (3') that  $d_{C^{*}}(v_{1})=2$, and  $ 2 \leqslant d_{C^{*}}(u_{i}) \leqslant 3$ for any $i \in \{1,p\}$. 
	If $d_{C^{*}}(u_{1})=d_{C^{*}}(u_{p}) =3$, then either $N_{C^{*}}(u_{1})=N_{C^{*}}(u_{p})$ or $N_{C^{*}}(u_{1}) \cup N_{C^{*}}(u_{p})=V(C^{*})$, 
	and so $d_{C^{*}}(u_{1},u_p)=3$ or $6$, 
	contrary to (2').
	Thus, without loss of generality, let $d_{C^{*}}(u_{p})=2$.
	In addition, we may let $w_1 \in N_{C^{*}}(v_1)$.

	By (3) of Lemma \ref{degree-3},  $w_2,\,w_6\not\in N_{C^{*}}(v_{1})$.
	Suppose that $N_{C^{*}}(v_{1})=\{w_{1},w_{3}\}$. Then we may let $\{w_{4},w_{6}\} \subseteq N_{C^{*}}(u_{1}) \subseteq \{w_{2},w_{4},w_6\}$,  and then $N_{C^{*}}(u_{p})=\{w_{2},w_{5}\}$. Clearly, we have $w_{1}w_{4} \in E(G)$ or $w_{3}w_{6} \in E(G)$. 
	It is easy to check that $v_{1}w_{1}w_{6}w_{5}w_4w_{3}v_{1}$ is a $6$-cycle, write $C'$, with a chord $w_{1}w_{4}$ or $w_{3}w_{6}$.
	Recalling that $H-v_1$ is connected, we get an optimal $(s-1)$-system of chorded cycles $\mathcal{C}'$  by replacing $C^{*}$ with $C'$. Now  $G-\mathcal{C}'$ has a component $(H-v_1)\sqcup w_2$ which contains a path of length $p=\ell(\mathcal{C})+1$, contrary to the assumption $(\ddag)$. 
	Similarly, if $N_{C^{*}}(v_{1})=\{w_{1},w_{5}\}$ then we get a contradiction. 
	Therefore, $N_{C^{*}}(v_{1})=\{w_{1},w_{4}\}$, and so $w_1w_4 \notin E(G)$ by (3) of Lemma \ref{degree-3}.

	Again by (3) of Lemma \ref{degree-3}, either $w_2w_5$ or $w_3w_6$
	is a chord of $C^*$. Then $N_{C^{*}}(u_{1})$ contains neither $\{w_{2},w_{5}\}$ nor
	$\{w_{3},w_{6}\}$. Thus $N_{C^{*}}(u_{1})$ intersects each of $\{w_{1},w_{4}\}$, $\{w_{2},w_{5}\}$ and
	$\{w_{3},w_{6}\}$ in at most one element. Since
	$C^{*}\sqcup u_1$ is triangle-free by Lemma \ref{degree-3}, it follows that
	$N_{C^{*}}(u_{1})\subseteq \{w_{1},w_{3},w_5\}$ or $\{w_{2},w_{4},w_6\}$.
	Similarly, $N_{C^{*}}(u_{p})\subseteq \{w_{1},w_{3},w_5\}$ or $\{w_{2},w_{4},w_6\}$.
	Since $d_{C^*}(u_1,u_p)=4$, without loss of generality, we let $N_{C^{*}}(u_{1})\subseteq \{w_{2},w_{4},w_6\}$ and
	$N_{C^{*}}(u_{p})\subseteq  \{w_{1},w_{3},w_5\}$.
	Then $\{w_2,w_6\}\subseteq N_{C^{*}}(u_{1})$, and    $ N_{C^{*}}(u_{p})=\{w_3,w_5\}$ as $d_{C^{*}}(u_{p})=2$.

	Suppose that $w_4\in N_{C^{*}}(u_{1})$. Then we have a chorded $6$-cycle
	$(C^*-w_4)\sqcup u_p$ with a chord $w_2w_5$ or $w_3w_6$.
	Recall that $u_p$ is not a cut vertex of $H$.
	Replacing $C^*$ with $(C^*-w_4)\sqcup u_p$, we get an optimal $(s-1)$-system of
	chorded cycles  say $\mathcal{C}_{u_p}$ from $\mathcal{C}$.
	Now $(H-u_p)\sqcup w_4$ is a component of $G-\mathcal{C}_{u_p}$ of maximal order,
	and $(H-u_p)\sqcup w_4$ contains a longer path $v_1w_4u_1\cdots u_{p-1}$ than $P$, 
	contrary to the assumption $(\ddag)$. Therefore, we have
	\begin{equation}\label{eq-222}
		N_{C^{*}}(v_1)=\{w_1,w_4\},\,\,\,\, N_{C^{*}}(u_{1})=\{w_2,w_6\},\,\,\,\,  N_{C^{*}}(u_{p})=\{w_3,w_5\}.
	\end{equation}
	
	\vskip 5pt

	Applying   Lemma \ref{no-Hamiltonian-path} to $H$, $P$ and $Q$, since $d_H(v_1)=2$, one of the following cases occurs:
	\begin{itemize}
		\item[(i)] $q\leqslant 2$ and $H=P\sqcup Q$;
		\item[(ii)] there exists $w\in V(H - (P \cup v_{1}))$ such that $d_{H}(w) \leqslant 2$, $u_{1}w, u_{p}w \notin E(H)$ and $H-w$ is connected.
	\end{itemize}
	
	Suppose first that (ii) occurs. Then   $d_{C^{*}}(w)=d_{C^{*}}(v_1)=2$  by  Claim \ref{comp-H-5} and (3') above.
	Recall that $ N_{C^{*}}(u_{1})=\{w_{2},w_{6}\}$ and $ N_{C^{*}}(u_{p})=\{w_{3},w_{5}\}$.
	It follows from Claim \ref{comp-H-5} that $N_{C^{*}}(w)=\{w_{1},w_{4}\}=N_{C^{*}}(v_{1})$, in particular,
	$d_{C^{*}}(v_{1},w)=2$. If $v_{1}w \in E(G)$, then  $C^{*}\sqcup wv_1$ contains a chorded $4$-cycle with a chord $v_{1}w$, contrary to Lemma \ref{C-mini}.
	Now  let   $v_{1}w \notin E(G)$. Then
	$\sum_{C\in\mathcal{C}\setminus C^{*}}d_C(v_{1},w)=d_{G}(v_{1},w)-d_{C^{*}}(v_{1},w)-d_{H}(v_{1},w) \geqslant 4s-2-4=4(s-2)+2$,
	and thus $d_{C'}(v_{1},w)  \geqslant 5$ for some $C'\in \mathcal{C}\setminus C^{*}$.
	By Lemma \ref{6-cycle},  $|V(C')|=6$, and either $d_{C'}(v_{1}) \geqslant 3$ or $d_{C'}(w)\geqslant 3$, a contradiction.

	In the following,  we suppose that $q\leqslant 2$ and $H=P\sqcup Q$, and arrive at a contradiction by investigating the role of vertex $u_2$.
	First, we claim that $H-u_2$ is connected. Suppose the contrary, noting that $P-u_2$ is connected as $d_P(u_1)=2$, then
	$N_Q(u_2)\ne \emptyset$, and $N_P(v)\subseteq\{u_2\}$ for any $v\in V(Q)$.
	If $q=1$ then 
	$d_H(v_1)=d_P(v_1)\leqslant 1$, which contradicts (2') above.
	Thus $q=2$. Since $d_H(v_1)=2$, we have $u_2\in N_P(v_1)$.
	Then we have a longer path $u_p\cdots u_2v_1v_2$   than $P$, contrary to the choice of $P$. 
	Therefore, $H-u_2$ is connected.
	
	Suppose  that $u_{2}v_{1} \in E(H)$.
	Then $C'=u_{1}w_{6}w_{5}w_{4}w_{3}w_{2}u_{1}$ is a chorded cycle with a chord $w_{2}w_{5}$ or $w_{3}w_{6}$.
	Replacing $C^{*}$ with $C'$, we get an optimal $(s-1)$-system  $\mathcal{C}'$ of  chorded cycles from
	$\mathcal{C}$. However,  $G-\mathcal{C}'$ has a component $(H-u_1)\sqcup w_1$, which contains a path of length $p=\ell(\mathcal{C})+1$,  contrary to the assumption $(\ddag)$. 
	
	The argument above implies that  $u_{2}v_{1} \notin E(H)$.
	In addition, if $q=2$ then   $u_{2}v_{2} \notin E(H)$, otherwise, we have a longer path $u_p\cdots u_2v_2v_1$   than $P$, a contradiction. 
	Then $d_{H}(u_{2})=d_{P}(u_{2})  \leqslant 3$, and so $d_{H}(u_{2},v_{1}) \leqslant 5$. 
	We have $\sum_{C\in \mathcal{C}}d_C(u_{2},v_{1})=d_{G}(u_{2},v_{1})-d_{H}(u_{2},v_{1}) \geqslant 4s-5=4(s-1)-1$.
	Suppose that there is  $C'\in \mathcal{C}$ such that $d_{C'}(u_{2},v_{1})\geqslant 5$. Then $|V(C')| =6$ by Lemma  \ref{6-cycle}, and so  $d_{C'}(u_{2})=3$ and $d_{C'}(v_{1})=2$ by (3') above.
	Note that $N_{C'}(u_{1}) \cap N_{C'}(u_{2}) =\emptyset$; otherwise, $C' \sqcup u_1u_2$ contains a chorded $5$-cycle, contrary to Lemma \ref{C-mini}.  Thus $N_{C'}(u_{1},v_{1}) \subseteq V(C') \setminus N_{C'}(u_{2})$, yielding $|N_{C'}(u_{1},v_{1})| \leqslant 3$, which contradicts (2') above.
	Therefore,    $d_{C}(u_{2},v_{1}) \in \{3,4\}$ for any chorded cycle $C\in \mathcal{C}$. 
	
	Again by Lemma \ref{C-mini}, we deduce  that $N_{C^{*}}(u_{1}) \cap N_{C^{*}}(u_{2}) =\emptyset$. Recall that $ N_{C^{*}}(u_{1})=\{w_{2},w_{6}\}$, and so $w_2,\,w_6 \notin N_{C^{*}}(u_2)$.
	Then, since  $d_{C}(u_{2},v_{1}) \geqslant 3$ and $ N_{C^{*}}(v_{1})=\{w_{1},w_{4}\}$, 
	either $w_3$ or $w_5$ is contained in $N_{C^*}(u_{2})$. 
	Let $w_j\in N_{C^*}(u_{2})$, where $j\in \{3,5\}$. Recall that $ N_{C^{*}}(u_{p})=\{w_3,w_5\}$.
	Noting that $H-u_1-u_2$ is connected, it follows that $H_j=(H-u_1u_2)\sqcup(w_3w_4w_5-w_j)$
	is connected. It is easily checked that
	$H_j$ contains a path $u_3u_4 \cdots u_pw_jw_4v_1$ of length $p=\ell(\mathcal{C})+1$.
	Suppose that $j=3$. 
	Then we have a $6$-cycle $C_3=u_{1}w_{6}w_{1}w_{2}w_{3}u_{2}u_{1}$  with a chord $u_{1}w_{2}$,
	and  an  optimal $(s-1)$-system  $\mathcal{C}_3$ of chorded cycles
	obtained from $\mathcal{C}$ by replacing $C^*$ with $C_3$.
	Now $H_3$ is a component of $G-\mathcal{C}_3$ of maximal order,
	and so $\ell(\mathcal{C}_3)\geqslant p=\ell(\mathcal{C})+1$, contrary to the assumption $(\ddag)$. 
	For   $j=5$, we get an  optimal $(s-1)$-system  $\mathcal{C}_5$ of chorded cycles
	from  $\mathcal{C}$ by replacing $C^*$ with a $6$-cycle $u_{1}w_{2}w_{1}w_{6}w_{5}u_{2}u_{1}$
	that has a chord $u_{1}w_{6}$. In this case, $H_5$ is a component of $G-\mathcal{C}_5$ of maximal order,
	which gives rise to a similar contradiction as above.
	This completes the proof of Theorem \ref{T4}.



\end{document}